\documentclass[10pt]{article}
\usepackage[shortlabels]{enumitem}
\usepackage{pb-diagram}
\usepackage{amsmath, amsfonts}
\usepackage{amssymb}
\usepackage{amscd}
\usepackage{color}

\usepackage{amsthm}
\usepackage{nccmath}
\usepackage[normalem]{ulem}
\usepackage{parskip}

\newtheorem{Df}{Definition}[section]
\newtheorem{Te}[Df]{Theorem}
\newtheorem{Po}[Df]{Proposition}
\newtheorem{Cr}[Df]{Corollary}
\newtheorem{Lm}[Df]{Lemma}
\newtheorem{Ca}[Df]{Claim}
\newtheorem{Cn}[Df]{Conjecture}
\newtheorem{Ex}[Df]{Example}
\newtheorem{Rm}[Df]{Remark}

\newcommand{\Bdf}{\begin{Df}}
\newcommand{\Edf}{\end{Df}}
\newcommand{\Bte}{\begin{Te}}
\newcommand{\Ete}{\end{Te}}
\newcommand{\Bpo}{\begin{Po}}
\newcommand{\Epo}{\end{Po}}
\newcommand{\Bcr}{\begin{Cr}}
\newcommand{\Ecr}{\end{Cr}}
\newcommand{\Blm}{\begin{Lm}}
\newcommand{\Elm}{\end{Lm}}
\newcommand{\Bca}{\begin{Ca}}
\newcommand{\Eca}{\end{Ca}}
\newcommand{\Bcn}{\begin{Cn}}
\newcommand{\Ecn}{\end{Cn}}
\newcommand{\Bex}{\begin{Ex}}
\newcommand{\Eex}{\end{Ex}}
\newcommand{\Brm}{\begin{Rm}}
\newcommand{\Erm}{\end{Rm}}
\newcommand{\Bdm}{{\it Proof.}\ }
\newcommand{\Edm}{\rule{2mm}{2mm}}

\newcommand{\pgq}{\geqslant}

\newcommand{\modb}{\text{\textsf{\upshape Mod}-}\ensuremath{B}}
\newcommand{\abimod}{\ensuremath{A}\text{-\textsf{\upshape Bimod}}}
\newcommand{\catvs}{\text{\textsf{\upshape Vect}\ensuremath{_k}}}

\newcommand{\Hom}{\mathop{\rm Hom}\nolimits}

\newcommand{\RHom}{\mathop{\rm RHom}\nolimits}

\newcommand{\HH}{\mathop{\rm HH}\nolimits}
\newcommand{\HK}{\mathop{\rm HK}\nolimits}

\begin{document}

\title{\bf{A cup-cap duality in Koszul calculus}}
\author{Roland Berger and Andrea Solotar
\thanks{\footnotesize The second mentioned author has been supported by the projects  UBACYT 20020130100533BA, PIP-CONICET 11220150100483CO and PICT 2015-0366. She is a research member of CONICET (Argentina) and a Senior Associate at ICTP.}}

\date{}

\maketitle

\begin{abstract}
We introduce a cup-cap duality in the Koszul calculus of $N$-homogeneous algebras defined by the first author~\cite{berger:Ncal}. As an application, we prove that the graded symmetry of the Koszul cap product is a consequence of the graded commutativity of the Koszul cup product. We propose a conceptual approach that may lead to a proof of the graded commutativity, based on derived categories in the framework of DG bimodules over DG algebras. Various enriched structures are developed in a weaker
 situation corresponding to $N>2$.
\end{abstract} 


\section{Introduction}

Just after the appearance of the Koszul calculus on quadratic algebras over a ground field~\cite{bls:kocal}, two motivating generalizations have allowed to enlarge the validity of this calculus to wider domains of applications.

The first one concerns a generalization to $N$-homogeneous algebras, that is, instead of the quadratic relations, homogeneous 
relations of degree $N>2$ over the ground field are considered~\cite{berger:Ncal}. This generalization is motivated by a great amount of works 
developed in the last two decades about $N$-Koszul algebras and various applications (see~\cite{berger:Ncal} for an extended bibliography).

In the second generalization, general quiver algebras with quadratic relations are considered, the quivers having finitely many vertices and arrows.
This generalization, applied to preprojective algebras, reveals some interesting interactions with representation theory and Poincar\'e Van den Bergh duality as presented in~\cite{bt}, with a particular interest on derived categories.

In the present paper, we have chosen to focus on $N$-homogeneous algebras over a ground field. A Koszul calculus for $N$-homogeneous relations on quivers is beyond our scope, even if such a calculus deserves to be defined and studied. Actually, in the present paper, we are essentially interested in a 
question that we detail below and which is already crucial in the one vertex case.

A classical result in Hochschild calculus asserts that the cup product is graded commutative on cohomology classes (Gerstenhaber's theorem~\cite{gerst:cohom}). By analogy with Hochschild calculus, a natural question in the Koszul calculus defined on $N$-homogeneous algebras is to know whether, on Koszul classes, the Koszul cup product $\underset{K}{\smile}$ is graded commutative. Similarly, one can ask whether the Koszul cap product $\underset{K}{\frown}$ is graded symmetric. These questions are still open when the algebra is quadratic, that is, when $N=2$, although the answers are positive if the quadratic algebra is Koszul~\cite[Subsection 3.4]{bls:kocal} and there are no examples up to date neither of non commutativity nor of non symmetry. When $N\geq 2$, we prove that the graded symmetry of the Koszul cap product is a consequence of the graded commutativity of the Koszul cup product. For that, we introduce a cup-cap duality in the Koszul calculus of $N$-homogeneous algebras.

For an $N$-homogeneous algebra $A$, our cup-cap duality consists in a duality linking a Koszul cup bracket $[\alpha , -]_{\underset{K}{\smile}}$ and a corresponding cap bracket $[\alpha , -]_{\underset{K}{\frown}}$. In the first argument of the brackets, $\alpha$ denotes any Koszul cohomology class with coefficients in $A$, while the second argument applies to any Koszul cohomology (respectively, homology) class with coefficients in an $A$-bimodule. The statement is the following and constitutes the main result of our paper (Theorem \ref{cupcaptheo} in Section 4). The various ingredients involved in this statement, in particular the duality isomorphism $\zeta$, will be thoroughly defined and studied in our paper.

\Bte \label{cupcapintro}
Let $A$ be an $N$-homogeneous algebra over a $k$-vector space $V$ and let $M$ be an $A$-bimodule. For any $q\geq p\geq 0$ and any $\alpha \in \HK^p(A)$, 
 there is a commutative diagram in the category of $k$-vector spaces
\begin{eqnarray} \label{cupcapduality6}
\HK_{q-p}(A,M)^{\ast}  \ \ \ \stackrel{[\alpha, -]_{\underset{K}{\frown}}^{\ast}}{\longrightarrow} & \HK_q(A,M)^{\ast} \nonumber  \\
\downarrow  \zeta_{q-p} \ \ \ \ \ \ \ \ \ \ \   &  \downarrow \zeta_q \\
\HK^{q-p}(A,M^{\ast}) \ \ \ \stackrel{- [\alpha , -]_{\underset{K}{\smile}}}{\longrightarrow} & \HK^q(A,M^{\ast}) \nonumber
\end{eqnarray}
where the vertical arrows are isomorphisms.
\Ete

Accordingly, it is immediate that if the bracket $[\alpha , -]_{\underset{K}{\smile}}$ acting on any $A$-bimodule is zero, then the bracket $[\alpha , -]_{\underset{K}{\frown}}$ acting on any $A$-bimodule is as well zero. We express this implication by saying that the graded symmetry of $\underset{K}{\frown}$ is a consequence of the graded commutativity of $\underset{K}{\smile}$.

We develop a derived category approach -that we hope will lead to a proof of the graded commutativity- inspired by the idea of ``enriched structures'' used 
in~\cite{bt}. Unfortunately, we have to add hypotheses concerning the existence of injective resolutions. However, we think that the  material so developed in the framework of derived categories is of interest in its own. 

The contents of the paper are as follows.
In Section \ref{recalling} we recall the basic definitions and results.
In Section \ref{Knatural}, we apply a general natural isomorphism to Koszul calculus. In Section~\ref{Kcupcap}, we present and prove Theorem \ref{cupcapintro} where the definition of the duality isomorphism $\zeta$ uses the isomorphism of Section \ref{Knatural}. As explained above, the graded symmetry of the Koszul cap product becomes a consequence of the graded commutativity of the Koszul cup product, and we also prove in Corollary~\ref{cupcapcons} that both properties hold if the algebra is $N$-Koszul.
We present in Section~\ref{examples} some examples of non $N$-Koszul algebras satisfying both properties.

In Section \ref{appli}, we verify the naturality in $A$ of the cup-cap duality isomorphism $\zeta$ specialized to $M=A$, and we combine it for $N=2$ with the Koszul duality obtained in~\cite{bls:kocal}. In Section \ref{Hcupcap}, we construct a cup-cap duality in Hochschild calculus for any associative algebra. Therefore our Hochschild cup-cap duality allows us to prove that the cap product on Hochschild classes is graded symmetric.

In Section \ref{enristruct}, we generalize to $N$-homogeneous algebras the enriched structure obtained on the Koszul bimodule complex $K(A)$ when $N=2$~\cite{bt}. In this generalization, 
there is a difficulty coming from the fact that the Koszul cup and cap products are not associative on cochains and chains when $N>2$~\cite{berger:Ncal}. So in this case 
we have to work with weak DG algebras and weak DG bimodules. It is noteworthy that the enriched structures obtained in~\cite{bt} when $N=2$ make sense in the 
weak situation $N>2$.

In Section \ref{towards}, we attack the question of the graded commutativity of the Koszul cup product. We base this section on a general 
result showing that the Koszul cohomology is isomorphic to a Hochschild hypercohomology (Proposition \ref{isohyper}). Then, using the enriched structure of $K(A)$ obtained in Section \ref{enristruct}, the three functors involved in Proposition \ref{isohyper} are enriched by suitable weak DG structures (Proposition \ref{enrichedrho}). We end Section \ref{towards} with an application (Theorem \ref{derivedrho}) which constitutes, in terms of derived categories, some evidence to state the graded commutativity when $N=2$.  

Throughout the article, we use notation and results from~\cite{berger:Ncal}. We fix a vector space $V$ over a field $k$ and an integer $N \geq 2$. The tensor algebra $T(V)=\bigoplus_{m\geq 0} V^{\otimes m}$ is graded by the \emph{weight} $m$. We also fix a subspace $R$ of $V^{\otimes N}$. The associative algebra $A=T(V)/(R)$ inherits the weight grading. The homogeneous component of weight $m$ of $A$ is denoted by $A_m$. In particular, $A_0=k$ and $A_m= V^{\otimes m}$ if $1 \leq m \leq N-1$. The graded algebra $A$ is called an \emph{N-homogeneous algebra}. If $N=2$, $A$ is called a quadratic algebra~\cite{bls:kocal}. In the next section, we recall some objects and results that are necessary in the sequel.

\setcounter{equation}{0}

\section{The Koszul calculus of $A$} \label{recalling}

\subsection{The Koszul bimodule complex of $A$} \label{Kbc}

Let $A=T(V)/(R)$ be an $N$-homogeneous algebra, where $R$ is a subspace of $V^{\otimes N}$. The fundamental object of the Koszul calculus of $A$ is the Koszul bimodule complex $K(A)$ that we recall below. For any $A$-bimodule $M$, the Koszul homology space $\HK_{\bullet}(A,M)$ is defined as the homology of the Koszul chain complex $(M\otimes W_{\nu(\bullet)}, b_K)$. In a similar way, the Koszul cohomology space $\HK^{\bullet}(A,M)$ is defined as the cohomology of the Koszul cochain complex $(\Hom(W_{\nu(\bullet)}, M), b_K)$. See~\cite[Subsection 2.2]{berger:Ncal} for the definitions of the Koszul differentials $b_K$ from the complex $K(A)$. 

For the definition of $K(A)$, we follow~\cite{berger:Ncal}. 
For any $p\geq N$, let $W_p$ be the subspace of $V^{\otimes p}$ defined by 
\[W_{p}=\bigcap_{i+N+j=p}V^{\otimes i}\otimes R\otimes V^{\otimes j}, \hbox{ where } i, j\ge 0,\]
\noindent while $W_0=k$ and $W_p=V^{\otimes p}$ if $1 \leq p \leq N-1$. It is convenient to use the following notation: an 
arbitrary element of $W_p$ will be denoted by a product $x_1 \ldots  x_p$, where $x_1, \ldots , x_p$ are in $V$. This notation should be thought of 
as a sum of such products. Moreover, regarding $W_p$ as a subspace of $V^{\otimes q}\otimes W_r \otimes V^{\otimes s}$ where $q+r+s=p$, the 
element $x_1 \ldots  x_p$ viewed in $V^{\otimes q}\otimes W_r \otimes V^{\otimes s}$ will be denoted in the same way, meaning that 
$x_{q+1} \ldots x_{q+r}$ is thought of as a sum belonging to $W_r$ and the other $x_i$'s are arbitrary elements in $V$. 

Let $\nu:\mathbb{N}\to \mathbb{N}$ be the map such that $\nu(p)=Np'$ if $p=2p'$ and $\nu(p)=Np'+1$ if $p=2p'+1$. The $A$-bimodule complex $K(A)$ is 
\begin{equation} \label{priddy}
\cdots \stackrel{d}{\longrightarrow} K_{p} \stackrel{d}{\longrightarrow} K_{p-1} \stackrel{d}{\longrightarrow} \cdots
\stackrel{d}{\longrightarrow} K_{1} \stackrel{d}{\longrightarrow} K_{0} \stackrel{}{\longrightarrow} 0\,,
\end{equation}
where $K_{p}=A\otimes W_{\nu(p)} \otimes A$. For any $a,a'\in A$ and $x_1\dots x_{\nu(p)} \in W_{\nu(p)}$, 
the differential $d$ is defined on $K_p$ as follows. If $p=2p'+1$, one has 
\[
d(a \otimes x_1 \ldots x_{Np'+1} \otimes a') =ax_1\otimes x_2 \ldots x_{Np'+1}\otimes a'- a \otimes x_1 \ldots x_{Np'}\otimes x_{Np'+1} a', \]
and if $p=2p'$, one has
\[
d(a \otimes x_1 \ldots x_{Np'} \otimes a') = \sum_{0\le i \le N-1}ax_1\ldots x_i\otimes x_{i+1}\ldots x_{i+Np'-N+1}\otimes x_{i+Np'-N+2}\ldots x_{Np'}a'.\]
The homology of $K(A)$ is equal to $A$ in degree $0$, and to $0$ in degree $1$. 

In case $N=2$, $W_p$ is a subspace of $V^{\otimes p} \subseteq A^{\otimes p}$ for all $p\geq 0$, and $K(A)$ is an $A$-bimodule subcomplex of the bar resolution $B(A)$ of $A$. However, if $N>2$, there is no natural inclusion of $K(A)$ in $B(A)$.

Koszul algebras for $N>2$ were defined in \cite{rb:NK} and the following
equivalent definition appeared in \cite{bm}. For $N=2$, this definition includes Priddy's definition~\cite{priddy} (see also~\cite{vdb:nch,bgs}).
\Bdf \label{defK}
An $N$-homogeneous algebra $A=T(V)/(R)$ is said to be $N$-Koszul if the homology of $K(A)$ is $0$ in any positive degree. 
\Edf

\subsection{Embedding $K(A)$ into a minimal resolution} \label{embedmin}

As recalled in~\cite[Subsection 2.6]{berger:Ncal}, in the category of graded $A$-bimodules, $A$ has a minimal projective resolution $P(A)$ whose component 
of homological degree $p$ has the form $A\otimes E_p \otimes A$, where $E_p$ is a graded space where the degree is still called weight. Moreover, the components of weight less than $\nu(p)$ in $E_p$ are zero while the component of weight $\nu(p)$ contains $W_{\nu(p)}$~\cite{rb:NK, hkl:NMMT}. So, for any $N\geq 2$, $K(A)$ can be viewed as a weight graded subcomplex of $P(A)$, and $A$ is $N$-Koszul if and only if $P(A)=K(A)$. Actually, $E_0=k$, $E_1=V$, $E_2=R$ and the component of weight $\nu(3)$ in $E_3$ is equal to $W_{\nu(3)}$. However, if $N>2$, the component of weight $\nu(4)$ in $E_4$ may strictly contain $W_{\nu(4)}$. It is well-known that, if $N=2$, the component of weight $p$ in $E_p$ is equal to $W_p$ for any $p$. See \cite[Subsections 2.7 and 2.8]{bgs} and \cite[Chapter 1, Proposition 3.1]{popo:quad}, see also \cite[Subsection 2.8]{bt} for a generalization to quadratic quiver algebras.

We use now the resolution $P(A)$ for computing the Hochschild (co)homology of $A$ and we obtain the following $N$-analogue of~\cite[Proposition 2.13]{bt} (limited to the one vertex case) by the same proof. 

\Bpo \label{comparison}
Let $A=T(V)/(R)$ be an $N$-homogeneous algebra. For any $A$-bimodule $M$, the inclusion of $K(A)$ into $P(A)$ induces a linear map $\HK_2(A,M) \rightarrow \HH _2(A,M)$ which is surjective, and a linear map $\HH ^2(A,M) \rightarrow \HK^2(A,M)$which is injective.
\Epo

As in~\cite[Corollary 2.14]{bt} (still limited to the one vertex case), we have a more precise result when $M=A$. In fact, for any $N$-homogeneous algebra $A$, we can define a coefficient weight grading in Koszul (co)homology with coefficients in $A$. It suffices to extend naturally the definition of this grading from the case $N=2$~\cite{bls:kocal, bt} to any $N>2$. The grading of $\HK _p(A)$ and $\HK ^p(A)$ by the coefficient weight $r$ is denoted by $\HK _p(A)_r$ and $\HK ^p(A)_r$. Unlike the Koszul differentials, the Hochschild differentials are not homogeneous for the coefficient weight, but only for the total weight. The grading of $\HH _p(A)$ and $\HH ^p(A)$ for the total weight $t$ is denoted by $\HH _p(A)_t$ and $\HH ^p(A)_t$. Then, specializing $M=A$ in Proposition \ref{comparison}, the linear map $\HK_2(A) \rightarrow \HH _2(A)$ is homogeneous from the coefficient weight $r$ to the total weight $r+N$, while the linear map $\HH ^2(A) \rightarrow \HK^2(A)$ is homogeneous from the total weight $r-N$ to the coefficient weight $r$.

\Bcr \label{2degree2}
Let $A=T(V)/(R)$ be an $N$-homogeneous algebra.
\begin{enumerate}[\itshape(i)]
\item The linear map $\HK_2(A)_r \rightarrow \HH _2(A)_{r+N}$ is an isomorphism if $r=0$ and $r=1$.

\item  Assume that $A$ is finite dimensional. Let $\max$ be the highest $m$ such that $A_m\neq
  0$.
  The linear map $\HH ^2(A)_{r-N} \rightarrow \HK^2(A)_r$ is an isomorphism if $r=\max$ and $r=\max-1$.
\end{enumerate}
\Ecr

\subsection{Products} \label{products}

We also recall definitions and properties of the Koszul cup and cap products as stated in~\cite[Sections 3 and 5]{berger:Ncal}.

\Bdf \label{defcup}
Let $A=T(V)/(R)$ be an $N$-homogeneous algebra. Let $P$ and $Q$ be $A$-bimodules. Given Koszul cochains $f:W_{\nu(p)}\rightarrow P$ and $g:W_{\nu(q)}\rightarrow Q$, we define the Koszul $(p+q)$-cochain $f\underset{K}{\smile} g : W_{\nu(p+q)} \rightarrow P\otimes_A Q$ by 
\\
1. if $p$ and $q$ are not both odd, so that $\nu(p+q)=\nu(p)+\nu(q)$, one has
\begin{equation*}
(f\underset{K}{\smile} g) (x_1 \ldots x_{\nu(p+q)}) = f(x_1 \ldots x_{\nu(p)})\otimes_A \, g(x_{\nu(p)+1} \ldots  x_{\nu(p)+\nu(q)}),
\end{equation*}
2. if $p$ and $q$ are both odd, so that $\nu(p+q)=\nu(p)+\nu(q)+N-2$, one has
\begin{eqnarray*}
    (f\underset{K}{\smile} g) (x_1 \ldots x_{\nu(p+q)})  =  -\sum_{0\leq i+j \leq N-2} x_1\ldots x_i f(x_{i+1} \ldots x_{i+\nu(p)})x_{i+\nu(p)+1}\ldots x_{\nu(p)+N-j-2} \\
  \otimes_A \ g(x_{\nu(p)+N-j-1} \ldots  x_{\nu(p)+\nu(q)+N-j-2})x_{\nu(p)+\nu(q)+N-j-1} \ldots  x_{\nu(p)+\nu(q)+N-2}.
\end{eqnarray*}
\Edf

\Bdf \label{defcap}
Let $A=T(V)/(R)$ be  an $N$-homogeneous algebra. Let $M$ and $P$ be $A$-bimodules. For any $p$-cochain $f:W_{\nu(p)}\rightarrow P$ and any $q$-chain $z=m \otimes x_1 \ldots x_{\nu(q)}$ 
in $M\otimes W_{\nu(q)}$, we define the $(q-p)$-chains $f \underset{K}{\frown} z$ and $z \underset{K}{\frown} f$ with coefficients in $P\otimes_A M$ and $M\otimes_A P$ respectively, as follows.
\\
1. If $p$ and $q-p$ are not both odd, so that $\nu(q-p)=\nu(q)-\nu(p)$, one has
\begin{eqnarray*}
f\underset{K}{\frown} z = (f(x_{\nu(q-p)+1} \ldots  x_{\nu(q)})\otimes_A m)\otimes \, x_1 \ldots x_{\nu(q-p)},\\
z\underset{K}{\frown} f = (-1)^{pq} (m \otimes_A f(x_1 \ldots  x_{\nu(p)}))\otimes \, x_{\nu(p)+1} \ldots x_{\nu(q)}. 
\end{eqnarray*}
2. If $p=2p'+1$ and $q=2q'$, so that $\nu(q-p)=\nu(q)-\nu(p)-N+2$, one has 
\begin{eqnarray*}
  f\underset{K}{\frown} z =  -\sum_{0\leq i+j \leq N-2} (x_{Nq'-Np'-N+i+2}\ldots x_{Nq'-Np'-j-1} f(x_{Nq'-Np'-j} \ldots x_{Nq'-j}) \\
  \otimes_A x_{Nq'-j+1}\ldots x_{Nq'}m x_1 \ldots x_i)\otimes \, x_{i+1} \ldots  x_{i+Nq'-Np'-N+1}.
\end{eqnarray*}
\begin{eqnarray*}
  z\underset{K}{\frown} f =  \sum_{0\leq i+j \leq N-2} (x_{Nq'-j+1}\ldots x_{Nq'}m x_1 \ldots x_i \otimes_A f(x_{i+1} \ldots x_{Np'+i+1})\\  x_{Np'+i+2}\ldots x_{Np'+N-j-1}) 
  \otimes \, x_{Np'+N-j}\ldots x_{Nq'-j}.
\end{eqnarray*}
The chain $f \underset{K}{\frown} z$ is called the left Koszul cap product of $f$ and $z$, while $z \underset{K}{\frown} f$ is called their right Koszul cap product.
\Edf

There is an important fact on these products that we recall now. For any Koszul cochains $f$, $g$, $h$ and any Koszul chain $z$, the {\em associativity relations}
\begin{eqnarray}
  (f\underset{K}{\smile}g)\underset{K}{\smile} h & = & f\underset{K}{\smile} (g\underset{K}{\smile} h), \label{assocrela1} \\
  f\underset{K}{\frown} (g\underset{K}{\frown} z) & = & (f\underset{K}{\smile}g)\underset{K}{\frown} z, \label{assocrela2} \\
  (z\underset{K}{\frown} g) \underset{K}{\frown} f & = & z \underset{K}{\frown} (g\underset{K}{\smile} f), \label{assocrela3} \\
  f\underset{K}{\frown} (z\underset{K}{\frown} g) & = & (f\underset{K}{\frown}z)\underset{K}{\frown} g, \label{assocrela4}
\end{eqnarray}
hold for $N=2$, but they are no longer true in general for $N>2$. However these associativity relations hold on Koszul classes for any $N\geq 2$. See~\cite[Sections 3 and 5]{berger:Ncal} for details.

For any Koszul cochains $f\in\Hom(W_{\nu(p)},P)$ and
$g\in\Hom(W_{\nu(q)},Q)$, we have the identity
\begin{equation}
b_K(f \underset{K}{\smile} g)=b_k(f) \underset{K}{\smile} g+(-1)^p f \underset{K}{\smile} b_K(g).\label{eq:dga}
\end{equation}
Therefore, specializing $P=Q=A$, we see that $(\Hom(W_{\nu(\bullet)},A), b_K, \underset{K}{\smile})$ is a DG algebra when $N=2$, but only a weak DG algebra when $N>2$. Here ``weak'' means that the product of the algebra is not necessarily associative. For any $N> 2$, following the case of quadratic quiver algebras~\cite{bt}, the weak DG algebra $\Hom(W_{\nu(\bullet)},A)$ is denoted by $\tilde{A}$. Note that $\tilde{A}$ is $k$-central. Moreover, formula
\eqref{eq:dga} shows that for any $A$-bimodule $M$, $\Hom(W_{\nu(\bullet)},M)$ is a weak DG $\tilde{A}$-bimodule for the Koszul cup actions, only weak since the associativity relation \eqref{assocrela1} does not necessarily hold if $N>2$. However ``weak'' can be removed when $N=2$.

Similarly, for any Koszul $p$-cochain $f$ and any Koszul $q$-chain $z$, the following formulas
\begin{equation} \label{Nkoldga}
b_K(f \underset{K}{\frown} z)= b_K(f) \underset{K}{\frown} z +(-1)^p f \underset{K}{\frown} b_K(z),
\end{equation}
\begin{equation} \label{Nkordga}
b_K(z \underset{K}{\frown} f)= b_K(z) \underset{K}{\frown} f +(-1)^q z \underset{K}{\frown} b_K(f).
\end{equation}
show that $M\otimes W_{\nu(\bullet)}$ is a weak DG $\tilde{A}$-bimodule for the Koszul cap actions -- again weak since the associativity relations 
\eqref{assocrela2}, \eqref{assocrela3} and \eqref{assocrela4} do not necessarily hold if $N>2$. However ``weak'' can be removed when $N=2$. 

For any $N\geq 2$, $H(\tilde{A})=\HK^{\bullet}(A)$ is a graded associative algebra for the Koszul cup product, so that $\HK^{\bullet}(A,M)$ and $\HK_{\bullet}(A,M)$ are graded $\HK^{\bullet}(A)$-bimodules for the Koszul cup and cap actions respectively.
 
It is worth noting that the formula of $\smile_K$ for monomial algebras coincides with the formula that is obtained by restricting to the Koszul complex the formula for the usual cup product obtained via Bardzell's resolution $P(A)$, which is endowed by a diagonal. More precisely,
since the total complex of the bicomplex $P(A)\otimes_AP(A)$  is also a projective resolution of $A\otimes_A A\cong A$ as $A$-bimodule, there is a comparison map 
 $\Delta: P(A) \to  P(A)\otimes_AP(A)$ lifting the identity and unique up to homotopy, that can be used to compute the cup product in Hochschild cohomology:
if $f$ is an $n$-cocycle and $g$ is an $m$-cocycle, then $f \smile g= \mu_A(f\otimes g) \Delta$. Restricting this map to the Koszul complex, 
we obtain the formula of Definition~\ref{defcup}. 

\setcounter{equation}{0}

\section{A natural isomorphism in Koszul calculus} \label{Knatural}

Let $A=T(V)/(R)$ be an $N$-homogeneous algebra, where $R$ is a subspace of $V^{\otimes N}$.
Let $f:W_{\nu(p)}\rightarrow A$ be a Koszul $p$-cocycle with coefficients in $A$. Denote by $\alpha$ the class of $f$ in $\HK^p(A)$. Recall that the maps
$$f \underset{K}{\frown} -, \ - \underset{K}{\frown}f : M\otimes W_{\nu(q)} \longrightarrow M\otimes W_{\nu(q-p)}$$
define a left and a right action on the complex $(M\otimes W_{\nu(\bullet)}, b_K)$, inducing a left and a right action of $\alpha$ on $\HK_{\bullet}(A,M)$ defined by the maps
$$\alpha \underset{K}{\frown} -, \ - \underset{K}{\frown}\alpha : \HK_q(A,M) \longrightarrow \HK_{q-p}(A,M),$$
so that $\HK_{\bullet}(A,M)$ is a graded $\HK^{\bullet}(A)$-bimodule.

Throughout this paper, the (graded) dual space of a (graded) vector space $E$ is denoted by $E^{\ast}$, and the transpose map of a (graded) linear map $u$ is denoted by $u^{\ast}$. On one hand, we see
$$\HK_{\bullet}(A, M)^{\ast}= \bigoplus_{q\geq 0} \HK_q(A, M)^{\ast}$$
as a natural graded $\HK^{\bullet}(A)$-bimodule. In fact, the transpose map
$$R(\alpha):=(\alpha \underset{K}{\frown} -)^{\ast} : \HK_{q-p}(A,M)^{\ast} \longrightarrow \HK_q(A,M)^{\ast}$$
defines a right action since it is easy to verify that $R(\beta \underset{K}{\smile} \alpha)= R(\alpha)\circ R(\beta)$ for any $\beta$ in $\HK^{\bullet}(A)$.
Similarly $L(\alpha)=(- \underset{K}{\frown}\alpha)^{\ast}$ defines a left action. The associativity formula $\alpha \underset{K}{\frown}(\gamma \underset{K}{\frown} \beta)= (\alpha \underset{K}{\frown}\gamma ) \underset{K}{\frown} \beta$ for $\gamma$ in $\HK_{\bullet}(A,M)$ shows that $R(\alpha)\circ L(\beta)= L(\beta) \circ R(\alpha)$, so that $\HK_{\bullet}(A,M)^{\ast}$ is a graded $\HK^{\bullet}(A)$-bimodule.

On the other hand, we consider the cochain complex
$$(M\otimes W_{\nu(\bullet)})^{\ast}=\bigoplus_{q\geq 0}(M\otimes W_{\nu(q)})^{\ast}$$
endowed with the differential
$$b_K^{\ast}: (M\otimes W_{\nu(q)})^{\ast} \rightarrow (M\otimes W_{\nu(q+1)})^{\ast}.$$
For brevity, we denote this complex by $C$. Recall that if $\phi\in (M\otimes W_{\nu(q)})^{\ast}$, one has
$$b_K^{\ast}(\phi) = -(-1)^q \phi \circ b_K.$$

Assume that $\phi$ is a $q$-cocycle of $C$, meaning that
$$\phi : M\otimes W_{\nu(\bullet)} \longrightarrow k(-q)$$
is a morphism of complexes, where $k(-q)$ is the trivial complex $k$ concentrated in homological degree $q$, so that the class $\overline{\phi}$ in $H^q(C)$ coincides with the homotopy class of the complex morphism $\phi$. Then
$$H_q(\phi) : H_q(M\otimes W_{\nu(\bullet)}) \longrightarrow H_q(k(-q)) \cong k$$
only depends on $\overline{\phi}$, allowing us to define a linear map
$$\xi_q: H^q(C) \longrightarrow \HK_{q}(A,M)^{\ast}$$
by $\xi_q (\overline{\phi})=H_q(\phi)$, and a graded linear map
$$\xi : H^{\bullet}(C)=H^{\bullet}((M\otimes W_{\nu(\bullet)})^{\ast}) \longrightarrow \HK_{\bullet}(A,M)^{\ast}.$$
It is a consequence of a general fact~\cite[Corollaire 1, p.84]{bou:alghom} that $\xi$ is an isomorphism. Moreover, $H^{\bullet}((M\otimes W_{\nu(\bullet)})^{\ast})$ and $\HK_{\bullet}(A,M)^{\ast}$ are natural in $M$ and $\xi$ is a natural isomorphism.

\Bpo  \label{Knaturaliso}
The space $H^{\bullet}(C)$ is naturally a graded $\HK^{\bullet}(A)$-bimodule and the map $\xi$ is a natural isomorphism of graded $\HK^{\bullet}(A)$-bimodules. 
\Epo
\Bdm
Let $f:W_{\nu(p)}\rightarrow A$ be a Koszul $p$-cocycle. The linear map
$$(f \underset{K}{\frown} -)^{\ast}: (M\otimes W_{\nu(q-p)})^{\ast} \rightarrow (M\otimes W_{\nu(q)})^{\ast}$$
sends $\phi$ to $\phi \cdot f:=\phi \circ (f \underset{K}{\frown} -)$.
Assume that $\phi$ is a $(q-p)$-cocycle of $C$, meaning that $\phi \circ b_k=0$. Since $f$ is a $p$-cocycle, we have
$$(f \underset{K}{\frown}-) \circ b_K = (-1)^p b_K \circ (f \underset{K}{\frown} -)$$
on $M\otimes W_{\nu(q+1)}$. Thus $\phi \circ (f \underset{K}{\frown} -) \circ b_K =0$, implying that $\phi \cdot f$ is a $q$-cocycle. Then it is easy to check that $\overline{\phi \cdot f}$ in $H^q(C)$ only depends on $\overline{\phi}$ in $H^{q-p}(C)$. So we have defined the linear map
$$H^q((f \underset{K}{\frown} -)^{\ast}): H^{q-p}(C) \rightarrow H^q(C)$$

The reader will easily verify that this map only depends on the class $\alpha$ of $f$ in $\HK^p(A)$, and therefore it can be written as 
$R'(\alpha): H^{q-p}(C) \rightarrow H^q(C)$ with $R'(\alpha)(\overline{\phi})=\overline{\phi \cdot f}$ also denoted by $\overline{\phi} \cdot \alpha$. Next it is straightforward to check that $(\overline{\phi} \cdot \alpha) \cdot \beta = \overline{\phi} \cdot (\alpha\underset{K}{\smile}\beta)$, showing that $H^{\bullet}(C)$ endowed with the action $R'$ is a graded right $\HK^{\bullet}(A)$-module.

Similarly, the linear map $L'(\alpha): H^{q-p}(C) \rightarrow H^q(C)$ is well-defined by $L'(\alpha)(\phi)=\overline{f \cdot \phi}$ where $f \cdot \phi=\phi \circ (- \underset{K}{\frown} f)$, and $H^{\bullet}(C)$ endowed with the action $L'$ is a graded left $\HK^{\bullet}(A)$-module. Moreover it is easy to verify that the actions $R'$ and $L'$ commute. Thus $H^{\bullet}(C)$ is a graded $\HK^{\bullet}(A)$-bimodule.

Let us prove now that $\xi$ is right $\HK^{\bullet}(A)$-linear for the above actions (the left linearity is similar). Keeping the above notation, we have 
$$\xi_q(\overline{\phi} \cdot \alpha)=H_q(\phi \cdot f)= H_q(\phi \circ (f \underset{K}{\frown} -))=H_{q-p}(\phi) \circ H_q(f \underset{K}{\frown} -)=H_{q-p}(\phi) \circ (\alpha \underset{K}{\frown} -)$$
which is equal to $(\alpha \underset{K}{\frown} -)^{\ast}(H_{q-p}(\phi))=\xi_{q-p}(\overline{\phi}) \cdot \alpha$. The naturality in $M$ is clear.
\Edm

\Brm \label{remark}
We can define the actions $\phi \cdot f:=\phi \circ (f \underset{K}{\frown} -)$ and $f \cdot \phi:=\phi \circ (- \underset{K}{\frown} f)$ for any cochain $\phi$ of $C$ and any Koszul cochain $f\in \Hom(W_{\nu(\bullet)}, A)$. Then $C$ is a weak DG $\tilde{A}$-bimodule (see the end of Subsection \ref{products} for the definition of the weak DG algebra $\tilde{A}$). Passing to cohomologies, we obtain the graded $\HK^{\bullet}(A)$-bimodule $H^{\bullet}(C)$ defined above. 
\Erm

\setcounter{equation}{0}

\section{Koszul cup-cap duality} \label{Kcupcap}

Let us fix an $N$-homogeneous algebra $A=T(V)/(R)$ and an $A$-bimodule $M$. As usual, $M^{\ast}=\Hom(M,k)$ is an $A$-bimodule for the actions defined by $(a.u.a')(m)=u(a'ma)$ for $a\in A$, $a'\in A$ and $m\in M$. For any $p\geq 0$, we define the linear map
$$\eta_p: \Hom(M\otimes W_{\nu(p)}, k) \longrightarrow \Hom(W_{\nu(p)}, M^{\ast})$$
by $\eta_p (\varphi)(x_1 \ldots x_{\nu(p)})(m)=\varphi(m\otimes x_1 \ldots x_{\nu(p)})$ with obvious notations. It is a standard fact that the maps $\eta_p$ are linear isomorphisms. Their direct sum 
$$\eta: (M\otimes W_{\nu(\bullet)})^{\ast} \longrightarrow \Hom(W_{\nu(\bullet)}, M^{\ast})$$
is a graded linear isomorphism. This isomorphism is natural in the $A$-bimodule $M$. 

\Bpo 
The map $\eta$ is an isomorphism from the complex $C=((M\otimes W_{\nu(\bullet)})^{\ast}, b_K^{\ast})$ to the complex $(\Hom(W_{\nu(\bullet)}, M^{\ast}), b_K)$, 
lifting the identity of $M^*$ and inducing a graded linear isomorphism $H(\eta):H^{\bullet}(C) \rightarrow \HK^{\bullet}(A,M^{\ast})$ natural in the $A$-bimodule $M$.
\Epo
\Bdm
It amounts to prove that the diagram
\begin{eqnarray} \label{etamorphism}
(M\otimes W_{\nu(p)})^{\ast}  \ \ \ \stackrel{\eta_p}{\longrightarrow} & \Hom(W_{\nu(p)},M^{\ast}) \nonumber  \\
\downarrow  b_K^{\ast} \ \ \ \ \ \ \ \ \ \ \   &  \downarrow b_K \\
(M\otimes W_{\nu(p+1)})^{\ast} \ \ \ \stackrel{\eta_{p+1}}{\longrightarrow} & \Hom(W_{\nu(p+1)}, M^{\ast}) \nonumber
\end{eqnarray}
is commutative. Let $\phi \in (M\otimes W_{\nu(p)})^{\ast}$ and $m\in M$. Assuming that $p$ is even, we have the equality 
$$b_K(\eta_p(\phi))(x_1 \ldots x_{\nu(p+1)})(m) =\eta_p(\phi)(x_1\ldots  x_{\nu(p)}) (x_{\nu(p+1)}m) - \eta_p(\phi)(x_2 \ldots x_{\nu(p+1)})(mx_1).$$
The right-hand side is equal to
$$\phi(x_{\nu(p+1)}m\otimes x_1\ldots  x_{\nu(p)})-\phi(mx_1\otimes x_2 \ldots x_{\nu(p+1)})=-\phi \circ b_K(m\otimes x_1 \ldots x_{\nu(p+1)})$$
while the left-hand side coincides with $b_K^{\ast}(\phi)(m\otimes x_1 \ldots x_{\nu(p+1)})$. Thus $b_K(\eta_p(\phi))=\eta_{p+1}(b_K^{\ast}(\phi))$ as expected. 

When $p$ is odd, we use a similar argument. We have the equality
$$b_K(\eta_p(\phi))(x_1 \ldots x_{\nu(p+1)})(m) = \sum_{0\leq i\leq N-1} \eta_p(\phi)(x_{i+1} \ldots x_{i+\nu(p)})(x_{i+\nu(p)+1}\ldots x_{\nu(p+1)}mx_1\ldots x_i)$$
whose right-hand side is equal to
$$\sum_{0\leq i\leq N-1} \phi(x_{i+\nu(p)+1}\ldots x_{\nu(p+1)}mx_1\ldots x_i \otimes x_{i+1} \ldots x_{i+\nu(p)})=\phi \circ b_K(m\otimes x_1 \ldots x_{\nu(p+1)}). \ \Edm$$

According to Remark \ref{remark} and the end of Subsection \ref{products}, the complexes $C=(M\otimes W_{\nu(\bullet)})^{\ast}$ and $\Hom(W_{\nu(\bullet)}, M^{\ast})$ are weak DG $\tilde{A}$-bimodules.

\Bpo 
The map $\eta$ is a natural isomorphism of weak DG $\tilde{A}$-bimodules, inducing a natural isomorphism $H(\eta)$ of graded $\HK^{\bullet}(A)$-bimodules.
\Epo
\Bdm
Let us prove that $\eta$ is right $\tilde{A}$-linear, meaning that, for any $q\geq p\geq 0$ and any Koszul $p$-cochain $f:W_{\nu(p)}\rightarrow A$, the diagram 
\begin{eqnarray} \label{cupcapduality1}
\Hom(M\otimes W_{\nu(q-p)},k)  \ \ \ \stackrel{(f \underset{K}{\frown} -)^{\ast}}{\longrightarrow} & \Hom(M\otimes W_{\nu(q)},k) \nonumber  \\
\downarrow  \eta_{q-p} \ \ \ \ \ \ \ \ \ \ \   &  \downarrow \eta_q \\
\Hom(W_{\nu(q-p)}, M^{\ast}) \ \ \ \stackrel{\pm (- \underset{K}{\smile} f)}{\longrightarrow} & \Hom(W_{\nu(q)}, M^{\ast}) \nonumber
\end{eqnarray}
is commutative, where $\pm=(-1)^{(q-p)p}$. Note that this sign is a Koszul sign.

The vertical arrows are isomorphisms. Define
$$\psi: \Hom(W_{\nu(q-p)}, M^{\ast}) \rightarrow \Hom(W_{\nu(q)}, M^{\ast})$$
by $\psi = \eta_q \circ (f \underset{K}{\frown} -)^{\ast} \circ (\eta_{q-p})^{-1}$.
Let $g\in \Hom(M\otimes W_{\nu(q-p)},k)$ and $h=(f \underset{K}{\frown} -)^{\ast}(g)$. For brevity, denote $\eta_{q-p}(g)$ by $g'$ and $\eta_q(h)$ by $h'$, so that $\psi(g')=h'$. One has $h= (-1)^{(q-p)p} g\circ (f \underset{K}{\frown} -)$ so that
$$h(z)=(-1)^{(q-p)p} g (f \underset{K}{\frown} z)$$
for every $z=m\otimes x_1 \ldots x_{\nu(q)} \in M\otimes W_{\nu(q)}$.

Assuming $p=2p'+1$ and $q=2q'$, we have $\nu(p)=Np'+1$ and $\nu(q)= Nq'$. On one hand, $h'(x_1 \ldots x_{Nq'})(m)=h(z)$ with
\begin{eqnarray*}
h(z) =  \sum_{0\leq i+j \leq N-2} g(x_{Nq'-Np'-N+i+2}\ldots x_{Nq'-Np'-j-1} f(x_{Nq'-Np'-j} \ldots x_{Nq'-j}) \\
  x_{Nq'-j+1}\ldots x_{Nq'}m x_1 \ldots x_i \otimes \, x_{i+1} \ldots  x_{i+Nq'-Np'-N+1}).
\end{eqnarray*}
On the other hand, $\nu(q-p)=Nq'-Np'-N+1$ and $g': W_{\nu(q-p)} \rightarrow M^{\ast}$ is defined by $g'(x_1 \ldots x_{Nq'-Np'-N+1})(m)=g(m\otimes x_1 \ldots x_{Nq'-Np'-N+1})$, so that
\begin{eqnarray*}
(g'\underset{K}{\smile} f) (x_1 \ldots x_{Nq'})  =  -\sum_{0\leq i+j \leq N-2} x_1\ldots x_i g'(x_{i+1} \ldots x_{i+Nq'-Np'-N+1}) \\ x_{i+Nq'-Np'-N+2}\ldots x_{Nq'-Np'-j-1} 
  f(x_{Nq'-Np'-j} \ldots  x_{Nq'-j})x_{Nq'-j+1} \ldots  x_{Nq'}
\end{eqnarray*}
provides 
\begin{eqnarray*}
  (g'\underset{K}{\smile} f) (x_1 \ldots x_{Nq'})(m)  =  -\sum_{0\leq i+j \leq N-2} g(x_{i+Nq'-Np'-N+2}\ldots x_{Nq'-Np'-j-1} \\ f(x_{Nq'-Np'-j} \ldots  x_{Nq'-j})
x_{Nq'-j+1} \ldots  x_{Nq'} mx_1\ldots x_i \otimes \, x_{i+1} \ldots x_{i+Nq'-Np'-N+1} ).
\end{eqnarray*}
We obtain $h'=-g'\underset{K}{\smile} f$ as expected in the case $p$ odd and $q$ even. 

In the other cases, we have $(-1)^{(q-p)p}=1$ and $h'(x_1 \ldots x_{\nu(q)})(m)=h(z)$ with
$$h(z)= g(f(x_{\nu(q-p)+1}\ldots x_{\nu(q)})m\otimes x_1\ldots x_{\nu(q-p)}),$$
while $(g'\underset{K}{\smile} f) (x_1 \ldots x_{\nu(q)})  = g'(x_{1} \ldots x_{\nu(q-p)}) f(x_{\nu(q-p)+1} \ldots  x_{\nu(q)})$ provides
$$(g'\underset{K}{\smile} f) (x_1 \ldots x_{\nu(q)})(m)  = g(f(x_{\nu(q-p)+1} \ldots  x_{\nu(q)})m \otimes x_{1} \ldots x_{\nu(q-p)}).$$
Therefore $h'=g'\underset{K}{\smile} f$ as expected.

Let us prove similarly the left linearity of $\eta$, that is, the commutativity of the diagram
\begin{eqnarray} \label{cupcapduality2}
\Hom(M\otimes W_{\nu(q-p)},k)  \ \ \ \stackrel{\pm (- \underset{K}{\frown} f)^{\ast}}{\longrightarrow} & \Hom(M\otimes W_{\nu(q)},k) \nonumber  \\
\downarrow  \eta_{q-p} \ \ \ \ \ \ \ \ \ \ \   &  \downarrow \eta_q \\
\Hom(W_{\nu(q-p)}, M^{\ast}) \ \ \ \stackrel{f\underset{K}{\smile} -}{\longrightarrow} & \Hom(W_{\nu(q)}, M^{\ast}) \nonumber
\end{eqnarray}
where $\pm=(-1)^{pq}$ (a Koszul sign). We use the same notations as above. Now $\psi$ is defined by
$$\psi = (-1)^{pq} \eta_q \circ (- \underset{K}{\frown} f)^{\ast} \circ (\eta_{q-p})^{-1}$$
and $h(z)=(-1)^{pq}(-1)^{(q-p)p} g (z \underset{K}{\frown} f) = (-1)^p g(z \underset{K}{\frown} f)$.

First assume that $p=2p'+1$ and $q=2q'$. We have 
\begin{eqnarray*}
h(z) =  - \sum_{0\leq i+j \leq N-2} g(x_{Nq'-j+1}\ldots x_{Nq'}m x_1 \ldots x_i f(x_{i+1} \ldots x_{Np'+i+1})\\  x_{Np'+i+2}\ldots x_{Np'+N-j-1} 
  \otimes \, x_{Np'+N-j}\ldots x_{Nq'-j}).
\end{eqnarray*}
Next from
\begin{eqnarray*}
    (f\underset{K}{\smile} g') (x_1 \ldots x_{Nq'})  =  -\sum_{0\leq i+j \leq N-2} x_1\ldots x_i f(x_{i+1} \ldots x_{i+Np'+1})\\ x_{i+Np'+2}\ldots x_{Np'+N-j-1}
\ g'(x_{Np'+N-j} \ldots  x_{Nq'-j})x_{Nq'-j+1} \ldots  x_{Nq'},
\end{eqnarray*}
we deduce that
\begin{eqnarray*}
    (f\underset{K}{\smile} g') (x_1 \ldots x_{Nq'})(m)  =  -\sum_{0\leq i+j \leq N-2} g(x_{Nq'-j+1} \ldots  x_{Nq'}
m x_1\ldots x_i \\ f(x_{i+1} \ldots x_{i+Np'+1})x_{i+Np'+2}\ldots x_{Np'+N-j-1} \otimes \, x_{Np'+N-j} \ldots  x_{Nq'-j}),
\end{eqnarray*}
so that $h'=f\underset{K}{\smile} g'$ if $p$ odd and $q$ even. 

In the other cases, $h(z)= (-1)^p g(mf(x_{1}\ldots x_{\nu(p)})\otimes x_{\nu(p)+1}\ldots x_{\nu(q)})$ while 
$$(f\underset{K}{\smile} g') (x_1 \ldots x_{\nu(q)})  = (-1)^{pq} f(x_{1} \ldots x_{\nu(p)}) g'(x_{\nu(p)+1} \ldots  x_{\nu(q)}),$$
$$(f\underset{K}{\smile} g') (x_1 \ldots x_{\nu(q)})(m)  = (-1)^{pq} g(mf(x_{1} \ldots  x_{\nu(p)}) \otimes x_{\nu(p)+1} \ldots x_{\nu(q)}).$$
We conclude that $h'=f\underset{K}{\smile} g'$.
\Edm

Recall that the Koszul cup bracket $[- ,- ]_{\underset{K}{\smile}}$ and the Koszul cap bracket $[- , -]_{\underset{K}{\frown}}$ are defined as the graded commutators of the actions $\underset{K}{\smile}$ and $\underset{K}{\frown}$~\cite[Sections 4 and 6]{berger:Ncal}. From the commutative diagrams (\ref{cupcapduality1}) and (\ref{cupcapduality2}), we deduce immediately the following.
\Bpo \label{bracketduality}
Let $A$ be an $N$-homogeneous algebra over a $k$-vector space $V$, and let $M$ be an $A$-bimodule. For any $q\geq p\geq 0$ and any Koszul $p$-cochain 
$f:W_{\nu(p)}\rightarrow A$, there is a commutative diagram
\begin{eqnarray} \label{cupcapduality3}
\Hom(M\otimes W_{\nu(q-p)},k)  \ \ \ \stackrel{[f, -]_{\underset{K}{\frown}}^{\ast}}{\longrightarrow} & \Hom(M\otimes W_{\nu(q)},k) \nonumber  \\
\downarrow  \eta_{q-p} \ \ \ \ \ \ \ \ \ \ \   &  \downarrow \eta_q \\
\Hom(W_{\nu(q-p)}, M^{\ast}) \ \ \ \stackrel{-[f,- ]_{\underset{K}{\smile}}}{\longrightarrow} & \Hom(W_{\nu(q)}, M^{\ast}) \nonumber
\end{eqnarray}
\Epo

The diagrams (\ref{cupcapduality1}) and (\ref{cupcapduality2}) pass to the respective cohomologies. Moreover we use the isomorphism
$$\xi : H^{\bullet}((M\otimes W_{\nu(\bullet)})^{\ast}) \longrightarrow \HK_{\bullet}(A,M)^{\ast}$$
defined in Section \ref{Knatural}, and we obtain a natural graded $\HK^{\bullet}(A)$-bimodule isomorphism
$$\zeta : \HK_{\bullet}(A,M)^{\ast} \longrightarrow \HK^{\bullet}(A,M^{\ast})$$
by $\zeta= H(\eta) \circ \xi^{-1}$. Consequently, for any $q\geq p\geq 0$ and any $\alpha \in \HK^p(A)$, we have the commutative diagrams
\begin{eqnarray} \label{cupcapduality4}
\HK_{q-p}(A,M)^{\ast}  \ \ \ \stackrel{(\alpha \underset{K}{\frown} -)^{\ast}}{\longrightarrow} & \HK_q(A,M)^{\ast} \nonumber  \\
\downarrow  \zeta_{q-p} \ \ \ \ \ \ \ \ \ \ \   &  \downarrow \zeta_q \\
\HK^{q-p}(A,M^{\ast}) \ \ \ \stackrel{\pm (- \underset{K}{\smile} \alpha)}{\longrightarrow} & \HK^q(A,M^{\ast}) \nonumber
\end{eqnarray}
and
\begin{eqnarray} \label{cupcapduality5}
\HK_{q-p}(A,M)^{\ast}  \ \ \ \stackrel{\pm (- \underset{K}{\frown} \alpha)^{\ast}}{\longrightarrow} & \HK_q(A,M)^{\ast} \nonumber  \\
\downarrow  \zeta_{q-p} \ \ \ \ \ \ \ \ \ \ \   &  \downarrow \zeta_q \\
\HK^{q-p}(A,M^{\ast}) \ \ \ \stackrel{\alpha \underset{K}{\smile}-}{\longrightarrow} & \HK^q(A,M^{\ast}) \nonumber
\end{eqnarray}
from which we deduce the following.
\Bte \label{cupcaptheo}
Let $A$ be an $N$-homogeneous algebra over a $k$-vector space $V$ and let $M$ be an $A$-bimodule. For any $q\geq p\geq 0$ and any $\alpha \in \HK^p(A)$, 
 there is a commutative diagram in the category of $k$-vector spaces
\begin{eqnarray} \label{cupcapduality6}
\HK_{q-p}(A,M)^{\ast}  \ \ \ \stackrel{[\alpha, -]_{\underset{K}{\frown}}^{\ast}}{\longrightarrow} & \HK_q(A,M)^{\ast} \nonumber  \\
\downarrow  \zeta_{q-p} \ \ \ \ \ \ \ \ \ \ \   &  \downarrow \zeta_q \\
\HK^{q-p}(A,M^{\ast}) \ \ \ \stackrel{- [\alpha , -]_{\underset{K}{\smile}}}{\longrightarrow} & \HK^q(A,M^{\ast}) \nonumber
\end{eqnarray}
where the vertical arrows are isomorphisms.
\Ete

\Bcr \label{cupcapcons}
For any $N$-homogeneous algebra $A$ over a $k$-vector space $V$, the assertion (i) implies the assertion (ii), where

(i) For any $A$-bimodule $M$, any $\alpha \in \HK^{\bullet}(A)$ and any $\beta \in \HK^{\bullet}(A,M)$, $[\alpha, \beta]_{\underset{K}{\smile}}=0$.

(ii) For any $A$-bimodule $M$, any $\alpha \in \HK^{\bullet}(A)$ and $\gamma \in \HK_{\bullet}(A, M)$, $[\alpha, \gamma]_{\underset{K}{\frown}} =0$.

\noindent
In particular, if $A$ is $N$-Koszul, then (i) and (ii) hold.
\Ecr
\Bdm
The implication $(i) \Rightarrow (ii)$ is an immediate consequence of the theorem. 
As proved by Herscovich~\cite[Theorem 4.5]{hers:using}, if $A$ is $N$-Koszul, then the graded algebras $\HH^{\bullet}(A)$ and $\HK^{\bullet}(A)$ are isomorphic. Moreover, the graded $\HH^{\bullet}(A)$-bimodule $\HH^{\bullet}(A,M)$ is isomorphic to the graded $\HK^{\bullet}(A)$-bimodule $\HK^{\bullet}(A,M)$. 
Indeed, this fact is proved by Herscovich for actions on homology but his proof extends similarly to actions on cohomology. Therefore, 
when $A$ is $N$-Koszul, (i) holds since it holds in Hochschild calculus by a classical result of Gerstenhaber~\cite{gerst:cohom}.
\Edm
\\

If we know that the $A$-bimodule $M$ is $\mathbb{Z}$-graded with respect to the weight grading of $A$, then $M^{\ast}$ and more generally the functor 
$\Hom(M, -)$ can be replaced by their graded versions, including a Koszul sign in the graded actions as usual. So we obtain graded versions w.r.t 
$\mathbb{Z}$-gradings of bimodules of the maps $\xi$, $\eta$ and $\zeta$, and graded versions of the commutative diagrams (\ref{cupcapduality1})-(\ref{cupcapduality6}). If moreover $M$ is locally finite -- that is, each component of $M$ is finite-dimensional -- and if $M$ is replaced by $M^{\ast}$, then $M^{\ast}$ can be replaced by $M$, which gives a true duality.  When the $A$-bimodules involved in Corollary \ref{cupcapcons} are all assumed to be locally finite $\mathbb{Z}$-graded, the assertions (i) and (ii) are equivalent.

\setcounter{equation}{0}

\section{Examples of graded commutativity} \label{examples}

The aim of this section is to present some examples of non $N$-Koszul algebras satisfying both assertions (i) and (ii) of Corollary \ref{cupcapcons}. 

\subsection{A general result on the cup bracket} \label{generalbracket}

Let us recall the following result from~\cite[Corollary 4.6]{berger:Ncal}.

\Bpo \label{kocupcommutation}
Let $A$ be an $N$-homogeneous algebra over a $k$-vector space $V$, and let $M$ be an $A$-bimodule. For any $\alpha \in \HK^p(A,M)$ with $p=0$ or $p=1$ and any $\beta \in \HK^q(A)$ with $q\geq 0$, one has
\begin{equation} \label{kocupbrazero}
[\alpha, \beta]_{\underset{K}{\smile}} =0.
\end{equation}
\Epo

\Bcr \label{generators}
If $A=T(V)/(R)$ is an $N$-homogeneous algebra such that the Koszul cohomology algebra $\HK^{\bullet}(A)$ is generated in degrees $0$ and $1$, then for any $\alpha \in \HK^{\bullet}(A)$ and any $\beta \in \HK^{\bullet}(A)$, one has $[\alpha, \beta]_{\underset{K}{\smile}} =0$.
\Ecr
\Bdm
Given $\alpha \in \HK^p(A)$ and $\beta \in \HK^q(A)$, by hypothesis we can write 
\[\alpha = \sum a_{j_0}^0{}_{\underset{K}{\smile}} a_{j_1}^1{}_{\underset{K}{\smile}} \dots {}_{\underset{K}{\smile}} a_{j_{p}}^1\]
where the sum is finite, $a_{j_0}^0$ belongs to $\HK^0(A)$ and all the $a_{j_l}^1$'s belong to $\HK^1(A)$. Thus the result 
follows from Proposition~\ref{kocupcommutation} and from the fact that the cup bracket $[ -, -]_{\underset{K}{\smile}}$ is a derivation in the first argument on cohomology classes.
\Edm
\\

As the next example shows, the condition assumed in Corollary \ref{generators} is sufficient but not necessary.  

Let $A=k\langle x, y \rangle /(yx, y^2-xy)$. Let us denote by $r_1$ the relation $yx=0$ and by $r_2$ the relation $y^2 -xy=0$. 
The Koszul complex $K(A)$ is
\begin{equation} 
0 {\longrightarrow} A\otimes W_3 \otimes A \stackrel{d_2}{\longrightarrow} A\otimes R \otimes A \stackrel{d_1}{\longrightarrow} A\otimes V \otimes A
\stackrel{d_0}{\longrightarrow} A\otimes A  \longrightarrow 0
\end{equation}
where $W_3$ has one generator that we will call $w$ and the differentials are:
\begin{itemize}
\item $ d_0(1\otimes v \otimes 1) =  v\otimes 1 - 1 \otimes v$  for $v\in V$, as usual, 
\item $d_1(1\otimes r_1 \otimes 1) = y\otimes x \otimes 1 + 1 \otimes y \otimes x$, and $d_1(1\otimes r_2 \otimes 1) = y\otimes y \otimes 1 + 
1 \otimes y \otimes y -  1 \otimes x \otimes y - x \otimes y \otimes 1$,
\item $d_2(1\otimes w \otimes 1) = y \otimes r_1 \otimes 1 - x \otimes r_1 \otimes 1 - 1 \otimes r_2 \otimes x $.
\end{itemize}

The algebra $A$ is not Koszul, since $K(A)$ is not exact in homological degree $2$. Indeed, the element 
$y\otimes r_2\otimes x - x\otimes r_2\otimes x - x^2\otimes r_1\otimes 1$ is in the kernel of $d_1$ but not in the image of $d_2$.

Computing the Koszul cohomology, we obtain the following results. 
\begin{itemize}
\item $\HK^0(A) = \HH^0(A) =Z(A) = \{ a_0+a_1x^2+a_2xy+ \sum_{i=3}^na_ix^i, \hbox{ with }n \geq 3 \hbox{ and } a_i\in k \}$.
\item $\HK^1(A)=\HH^1(A)$ is generated by the classes of the elements $x^* \otimes x + y^* \otimes y$ and $x^* \otimes x^i$, for all $i\ge 2$.  
\item $\HK^2(A)$ is $1$-dimensional, generated by the class of $r_1^*\otimes y$. Note that $\HH^2(A)$ is also $1$-dimensional.
\item $\HK^3(A)$ is $1$-dimensional, generated by the class of $w^*\otimes xy$. We note that $\HH^3(A)$ is $2$-dimensional.
\item  Since the Koszul complex is zero in degrees greater or equal to $4$, we have that $\HK^i(A)=0$ for $i\ge 4$.
\end{itemize}

Then the computation of the Koszul cup products of these classes in positive degrees show that they are all zero. In particular, the algebra $\HK^{\bullet}(A)$ cannot be generated in degrees $0$ and $1$. Moreover, this algebra is graded commutative by Proposition \ref{kocupcommutation} with $M=A$.

\subsection{When the length of K(A) is at most 2}

If the Koszul bimodule complex $K(A)=A \otimes W_{\nu(\bullet)} \otimes A$ has length at most $2$, that is, if $W_{\nu(p)}=0$ whenever $p\geq 3$, 
then $\HK^p(A,M)=0$ for all $p\geq 3$ and therefore the cup actions of $\HK^q(A)$ on $\HK^p(A,M)$ are zero whenever $p+q \geq 3$. So Proposition \ref{kocupcommutation} and Corollary \ref{cupcapcons} imply the following.

\Bpo \label{lengthtwo}
Let $A$ be an $N$-homogeneous algebra over a $k$-vector space $V$. If the length of the Koszul bimodule complex $K(A)$ is at most $2$, 
then both assertions (i) and (ii) of Corollary \ref{cupcapcons} hold.
\Epo
\Bdm
It suffices to prove the equality (\ref{kocupbrazero}) when $p=2$ and $q=0$. But $\HK^0(A)=Z(A)$ the center of $A$, and $Z(A)$ acts symmetrically on $\HK^{\bullet}(A,M)$.
\Edm
\\

We apply now this proposition to the class of $N$-homogeneous algebras defined by a single monomial relation, as introduced in~\cite[Proposition 4.2]{rb:gera}.

\Bpo \label{onemon}
Let $V$ be a non-zero $k$-vector space of finite dimension $n$. Fix a basis $(x_1, \ldots, x_n)$ of $V$. 
Let $f=x_{i_1} \ldots x_{i_N}$ be a monomial of degree $N\geq 2$, and let $R$ be the subspace of $V^{\otimes N}$ 
generated by $f$. Then the $N$-homogeneous algebra $A=T(V)/(R)$ satisfies both assertions (i) and (ii) of Corollary \ref{cupcapcons}.
\Epo
\Bdm
It is proved in~\cite[Lemma 4.3]{rb:gera} that $W_p=0$ for every $p\geq N+1$, except in the case $f=x_i^N$ for some $i$, $1\leq i \leq n$. 
If $f\neq x_i^N$, then $K(A)$ has length $2$ and we conclude by Proposition \ref{lengthtwo}. If $f=x_i^N$, then $A$ is $N$-Koszul 
(see Proposition below), and 
we conclude using Corollary \ref{cupcapcons}.
\Edm
\\

Within the class of $N$-homogeneous algebras $A$ defined by a single monomial relation, the non $N$-Koszul algebras are characterized as follows~\cite[Proposition 4.2]{rb:gera}.

\Bpo \label{nonkonemon}
Let $V$ be a non-zero $k$-vector space of finite dimension $n$. Fix a basis $(x_1, \ldots, x_n)$ of $V$. 
Let $f=x_{i_1} \ldots x_{i_N}$ be a monomial of degree $N\geq 2$, and let $R$ be the subspace of $V^{\otimes N}$ 
generated by $f$. Then $A=T(V)/(R)$ is not $N$-Koszul if and only if there exists $m$ in $\{2,\ldots , N-1\}$ such that 
$$f=(x_{i_1} \ldots x_{i_m})^q\, x_{i_1} \ldots x_{i_r},$$
where $N=mq+r$ with $0\le r < m$, and where $i_1,\ldots ,i_m$ are not all equal.
\Epo

For example, if $n=2$ and $(x,y)$ is a basis of $V$, then $A=k\langle x,y \rangle/(xyx)$ is not $3$-Koszul. Note that the algebras $A$ considered
in Proposition \ref{nonkonemon} are always Koszul when $N=2$.

In Subsection \ref{embedmin}, we have recalled how the Koszul complex $K(A)$ can be embedded into a minimal resolution $P(A)$. We want to examine this embedding in the situation of Proposition \ref{nonkonemon} when $A$ is not $N$-Koszul. In this case, the Koszul complex is zero in degrees greater or equal to $3$.
This is not the case for the minimal resolution $P(A)$. For monomial algebras, $P(A)$ is isomorphic to Bardzell's resolution \cite{Ba}. Recall that Bardzell's resolution can be written down 
\begin{equation} 
\cdots {\longrightarrow} A\otimes R_3 \otimes A \stackrel{d_2}{\longrightarrow} A\otimes R_2 \otimes A \stackrel{d_1}{\longrightarrow} A\otimes V \otimes A
\stackrel{d_0}{\longrightarrow} A\otimes A  \longrightarrow 0
\end{equation}
where $R_2$ is a generating set of the monomial relations and for all $i\ge 3$, $R_i$ is the set of $(i-1)$-ambiguities. In the example $A=k\langle x,y \rangle/(xyx)$, one has $R_i=\{(xy)^{i-1}x \}$ for all $i \geq 2$. In general, for $A$ monomial, $W_3$ is the subset of $R_3$ of diagonal elements \cite{Ba}.

\setcounter{equation}{0}

\section{Application to functoriality and to Koszul duality} \label{appli}

Denote by $\mathcal{C}$ the generalized Manin category of $N$-homogeneous $k$-algebras~\cite{bdvw:homog, berger:Ncal} and by $\mathcal{E}$ the category of graded $k$-vector spaces. Recall that in $\mathcal{C}$, the objects are the $N$-homogeneous algebras and the morphisms are the morphisms of graded algebras. The $A$-bimodule $A^{\ast}=\Hom(A,k)$ is defined by the actions $(a.u.a')(x)=u(a'xa)$ for any linear map $u:A\rightarrow k$, and $x$, $a$, $a'$ in $A$. We have the following result~\cite[Proposition 2.3]{berger:Ncal}. Notice that in this statement, $A^{\ast}$ can be replaced by the graded dual -- including a Koszul sign in the left action -- of the $A$-bimodule $A$ endowed with the weight grading.

\Bpo  \label{functorial}
The rules $A\mapsto \HK_{\bullet}(A)$ and $A\mapsto \HK^{\bullet}(A,A^{\ast})$ respectively define a covariant functor $X$ and a contravariant functor $Y$ from $\mathcal{C}$ to $\mathcal{E}$. 
\Epo

Let us specialize the isomorphism $\zeta$ defined just before Theorem \ref{cupcaptheo} to $M=A$. We obtain a graded $\HK^{\bullet}(A)$-bimodule isomorphism
$$\zeta : \HK_{\bullet}(A)^{\ast} \longrightarrow \HK^{\bullet}(A,A^{\ast}),$$
and it is easy to check that it is natural in $A$. In other words, if we denote by $L$ the endofunctor $E \mapsto E^{\ast}$ (graded dual) of $\mathcal{E}$, then $\zeta$ defines an isomorphism of functors $L \circ X \cong Y$ from $\mathcal{C}$ to $\mathcal{E}$. When $A^{\ast}$ is replaced by the graded dual, one has a graded version of $\zeta$ satisfying an analogous natural property.

We want to combine the graded version of $\zeta$ with the Koszul duality developed for the quadratic Koszul calculus in~\cite[Section 8]{bls:kocal}. In particular, $A^{\ast}$ denotes the graded dual of the weight-graded $A$-bimodule $A$, as in~\cite[Section 8]{bls:kocal}. Using the same notations and hypotheses as in~\cite[Section 8]{bls:kocal}, we suppose that $V$ is finite-dimensional, $N=2$ and $A=T(V)/(R)$ is a quadratic algebra. We have the endofunctor $D: A \mapsto A^!$ of the Manin category $\mathcal{C}$ of quadratic $k$-algebras over finite-dimensional vector spaces. Recall that we have an isomorphism
$$\theta :  \HK_{\bullet}(A) \cong \tilde{\HK}^{\bullet}(A^!, A^{!\ast})$$
from the $(\HK^{\bullet}(A),\underset{K}{\smile})$-bimodule $\HK_{\bullet}(A)$ with actions $\underset{K}{\frown}$, 
$\mathbb{N}\times \mathbb{N}$-graded by the biweight, to the $(\tilde{\HK}^{\bullet}(A^!),\tilde{\underset{K}{\smile}})$-bimodule
 $\tilde{\HK}^{\bullet}(A^!, A^{!\ast})$  with actions $\tilde{\underset{K}{\smile}}$, $\mathbb{N}\times \mathbb{N}$-graded by the inverse biweight. This statement uses the fact that the bigraded algebras $(\HK^{\bullet}(A),\underset{K}{\smile})$ and $(\tilde{\HK}^{\bullet}(A^!),\tilde{\underset{K}{\smile}})$ are isomorphic. See~\cite[Section 8]{bls:kocal} for details and proofs.
For any $p\geq 0$ and $m\geq 0$, one has a linear isomorphism
\begin{equation} \label{linearisohkduality}
\theta _{p,m} : \HK_p(A)_m \cong \tilde{\HK}^m(A^!, A^{!\ast})_p.
\end{equation} 
As noted in~\cite[Remark 8.10]{bls:kocal}, $\theta$ defines an isomorphism of functors $X \cong \tilde{Y} \circ D$ from $\mathcal{C}$ to $\mathcal{E}$, where $\tilde{Y}$ is the contravariant functor $ A\mapsto \tilde{\HK}^{\bullet}(A, A^{\ast})$.

Under the above assumptions, the maps $\xi : H^{\bullet}((A\otimes W_{\nu(\bullet)})^{\ast}) \longrightarrow \HK_{\bullet}(A)^{\ast}$ and $\eta: (A\otimes W_{\nu(\bullet)})^{\ast} \longrightarrow \Hom(W_{\nu(\bullet)}, A^{\ast})$ are homogeneous w.r.t. the biweight $(p,m)$. Therefore the same is true for the map 
$\zeta : \HK_{\bullet}(A)^{\ast} \longrightarrow \HK^{\bullet}(A,A^{\ast})$. Then $\zeta$ is the direct sum of linear isomorphisms
$$\zeta_{p,m} : \HK_p(A)_m^{\ast} \longrightarrow \HK^p(A,A^{\ast})_m.$$
Since the vector spaces involved in $\zeta_{p,m}$ are finite-dimensional, we can consider the graded $\HK^{\bullet}(A)$-bimodule isomorphism
$$\zeta^{\ast} : \HK^{\bullet}(A,A^{\ast})^{\ast} \longrightarrow \HK_{\bullet}(A),$$
direct sum of the linear isomorphisms
$$\zeta^{\ast}_{p,m} : \HK^p(A,A^{\ast})_m^{\ast} \longrightarrow  \HK_p(A)_m.$$

Combining cup-cap duality and Koszul duality, we define an isomorphism
$$\Theta := \theta \circ \zeta^{\ast} : \HK^{\bullet}(A,A^{\ast})^{\ast} \longrightarrow \tilde{\HK}^{\bullet}(A^!, A^{!\ast})$$
from the $(\HK^{\bullet}(A),\underset{K}{\smile})$-bimodule $\HK^{\bullet}(A,A^{\ast})^{\ast}$ with transpose actions of $\underset{K}{\smile}$, 
$\mathbb{N}\times \mathbb{N}$-graded by the biweight, to the $(\tilde{\HK}^{\bullet}(A^!),\tilde{\underset{K}{\smile}})$-bimodule
 $\tilde{\HK}^{\bullet}(A^!, A^{!\ast})$  with actions $\tilde{\underset{K}{\smile}}$, $\mathbb{N}\times \mathbb{N}$-graded by the inverse biweight.
For any $p\geq 0$ and $m\geq 0$, one has a linear isomorphism
$$\Theta_{p,m} : \HK^p(A,A^{\ast})_m^{\ast} \longrightarrow \tilde{\HK}^m(A^!, A^{!\ast})_p.$$
In conclusion, $\Theta$ defines an isomorphism of functors $L \circ Y \cong \tilde{Y} \circ D$ from $\mathcal{C}$ to $\mathcal{E}$.

Similarly, define an isomorphism
$$\Theta':= \tilde{\zeta}^{-1} \circ \theta \circ D : \HK_{\bullet}(A^!) \longrightarrow \tilde{\HK}_{\bullet}(A)^{\ast}$$
from the $(\HK^{\bullet}(A^!),\underset{K}{\smile})$-bimodule $\HK_{\bullet}(A^!)$ to the $(\tilde{\HK}^{\bullet}(A),\tilde{\underset{K}{\smile}})$-bimodule $\tilde{\HK}_{\bullet}(A)^{\ast}$. For any $p\geq 0$ and $m\geq 0$, one has a linear isomorphism
$$\Theta'_{p,m} : \HK_p(A^!)_m \longrightarrow \tilde{\HK}_m(A)^{\ast}_p.$$
 Then $\Theta'$ defines an isomorphism of functors $X \circ D \cong L \circ \tilde{X}$ from $\mathcal{C}$ to $\mathcal{E}$, where $\tilde{X}: A \mapsto \tilde{\HK}_{\bullet}(A)$.

\setcounter{equation}{0}

\section{A cup-cap duality in Hochschild calculus} \label{Hcupcap}

It is known by the experts that the cap product is graded symmetric on Hochschild homology classes. Because of lack of a suitable reference, we include 
here a proof of this result, that we present as a consequence of a Hochschild cup-cap duality. In fact the cup-cap duality in Koszul calculus makes sense in a similar manner in Hochschild calculus for any non necessarily graded algebra.

Let $A$ be a unital associative $k$-algebra. In Hochschild calculus, the definition of the cup product $\smile$ and the definition of the left and right cap product $\frown$ coincide with the definitions in Koszul calculus when $N=2$ and when the elements $x_i$ of $V$ involved in the spaces $W_p$ are chosen 
arbitrarily in $A$. The defining formulas of the cup and cap products for a Hochschild $p$-cochain $f : A^{\otimes p} \rightarrow P$, a Hochschild $q$-cochain $g : A^{\otimes q} \rightarrow Q$ and a Hochschild $q$-chain $z=m \otimes a_1 \ldots a_q \in M\otimes A^{\otimes q}$ with coefficients in $A$-bimodules $P$, $Q$ and $M$ respectively, are the following.
\begin{equation} \label{cup}
(f\smile g) (a_1 \ldots a_{p+q}) = (-1)^{pq} f(a_1 \ldots a_p)\otimes_A \, g(a_{p+1} \ldots  a_{p+q}),
\end{equation}
\begin{equation} \label{lcap}
f\frown z = (-1)^{(q-p)p} (f(a_{q-p+1} \ldots a_q)\otimes_A m) \otimes a_1 \ldots  a_{q-p},
\end{equation}
\begin{equation} \label{rcap}
  z\frown f = (-1)^{pq} (m\otimes_A f(a_1 \ldots a_p)) \otimes a_{p+1} \ldots  a_q.
\end{equation}

Denote by $b$ the Hochschild differentials. For any $A$-bimodule $M$, we still have an isomorphism of complexes
$$\eta^H: (\Hom(M\otimes A^{\otimes \bullet}, k), b^{\ast}) \rightarrow (\Hom(A^{\otimes \bullet}, M^{\ast}), b)$$
defined in the same way than $\eta$, that is,
$$\eta_q^H (\varphi)(a_1 \ldots a_q)(m)=\varphi (m\otimes a_1 \ldots a_q)$$
where $\varphi: M\otimes A^{\otimes q} \rightarrow k$, the $a_i$'s are in $A$ and $m\in M$.
Then it is easy to prove that, for $q\geq p\geq 0$ and a Hochschild $p$-cochain $f:A^{\otimes p} \rightarrow A$, the diagrams 
\begin{eqnarray} \label{Hcupcapduality1}
\Hom(M\otimes A^{\otimes q-p},k)  \ \ \ \stackrel{(f \frown -)^{\ast}}{\longrightarrow} & \Hom(M\otimes A^{\otimes q},k) \nonumber  \\
\downarrow  \eta^H_{q-p} \ \ \ \ \ \ \ \ \ \ \   &  \downarrow \eta^H_q \\
\Hom(A^{\otimes q-p}, M^{\ast}) \ \ \ \stackrel{\pm (- \smile f)}{\longrightarrow} & \Hom(A^{\otimes q}, M^{\ast}) \nonumber
\end{eqnarray}
\begin{eqnarray} \label{Hcupcapduality2}
\Hom(M\otimes A^{\otimes q-p},k)  \ \ \ \stackrel{\pm (- \frown f)^{\ast}}{\longrightarrow} & \Hom(M\otimes A^{\otimes q},k) \nonumber  \\
\downarrow  \eta^H_{q-p} \ \ \ \ \ \ \ \ \ \ \   &  \downarrow \eta^H_q \\
\Hom(A^{\otimes q-p}, M^{\ast}) \ \ \ \stackrel{f\smile -}{\longrightarrow} & \Hom(A^{\otimes q}, M^{\ast}) \nonumber
\end{eqnarray}
commute, where $\pm=(-1)^{(q-p)p}$ in (\ref{Hcupcapduality1}) and $\pm=(-1)^{pq}$ in (\ref{Hcupcapduality2}). In fact, following the proof of the commutativity of 
the diagrams (\ref{cupcapduality1}) and (\ref{cupcapduality2}), it suffices to assume that $N=2$ in this proof (no cases to distinguish) and that the elements 
$x_i$ are chosen arbitrarily in $A$. We leave the details to the reader.

We also have a natural graded $\HH^{\bullet}(A)$-bimodule isomorphism
$$\xi^H : H^{\bullet}((M\otimes A^{\bullet})^{\ast}) \longrightarrow \HH_{\bullet}(A,M)^{\ast}$$
defined as in Section \ref{Knatural} and we introduce a natural graded $\HH^{\bullet}(A)$-bimodule isomorphism
$$\zeta^H : \HH_{\bullet}(A,M)^{\ast} \longrightarrow \HH^{\bullet}(A,M^{\ast})$$
by $\zeta^H= H(\eta^H) \circ (\xi^H)^{-1}$.
The commutative diagrams (\ref{Hcupcapduality1}) and (\ref{Hcupcapduality2}) pass to cohomology, and for any $q\geq p\geq 0$, $\alpha \in \HH^p(A)$, we have the commutative diagram
\begin{eqnarray} \label{Hcupcapduality3}
\HH_{q-p}(A,M)^{\ast}  \ \ \ \stackrel{[\alpha, -]_{\frown}^{\ast}}{\longrightarrow} & \HH_q(A,M)^{\ast} \nonumber  \\
\downarrow  \zeta^H_{q-p} \ \ \ \ \ \ \ \ \ \ \   &  \downarrow \zeta^H_q \\
\HH^{q-p}(A,M^{\ast}) \ \ \ \stackrel{- [\alpha , -]_{\smile}}{\longrightarrow} & \HH^q(A,M^{\ast}) \nonumber
\end{eqnarray}

We know from Gerstenhaber~\cite{gerst:cohom} that $[\alpha, \beta]_{\smile}=0$ for any $A$-bimodule $M$, any $\alpha \in \HH^{\bullet}(A)$ and any $\beta \in \HH^{\bullet}(A,M)$. From the commutative diagram (\ref{Hcupcapduality3}), we thus obtain the expected result.

\Bte \label{Hcupcapcons}
Let $A$ be a unital associative algebra. Then the cap product is graded symmetric on Hochschild homology classes, that is, for any $A$-bimodule $M$, any $\alpha \in \HH^{\bullet}(A)$ and any $\gamma \in \HH_{\bullet}(A, M)$, we have $[\alpha, \gamma]_{\frown} =0$.
\noindent
\Ete

The application to functoriality developed in Section \ref{appli} follows along the same lines. Denote by $\mathcal{A}$ the category of unital associative $k$-algebras. Recall that $\mathcal{E}$ is the category of graded $k$-vector spaces. We have functors $X^H:A\mapsto \HH_{\bullet}(A)$ and $Y^H:A\mapsto \HH^{\bullet}(A,A^{\ast})$ from $\mathcal{A}$ to $\mathcal{E}$~\cite{loday:cychom}. Then the graded $\HH^{\bullet}(A)$-bimodule isomorphism
$$\zeta^H : \HH_{\bullet}(A)^{\ast} \longrightarrow \HH^{\bullet}(A,A^{\ast})$$
defines an isomorphism of functors $L \circ X^H \cong Y^H$ from $\mathcal{A}$ to $\mathcal{E}$.

\setcounter{equation}{0}

\section{An enriched structure on $K(A)$} \label{enristruct}

The motivation of this section comes from a general strategy, developed in the next section, for proving the graded commutativity of the Koszul cup product, that is, property (i) in Corollary \ref{cupcapcons}. This general strategy is inspired by a similar strategy used in~\cite{bt}, where a Koszul complex 
Calabi-Yau property is introduced for any quadratic quiver algebra $A$, this property implying a Poincar\'e Van den Bergh duality for the Koszul homology/cohomology of $A$ with coefficients in any $A$-bimodule $M$. In order to obtain a more precise duality expressed as a cap product by a fundamental class, a stronger version of the Koszul complex Calabi-Yau property is defined in~\cite[Section 5]{bt}. The stronger definition consists in enriching the DG $\tilde{A}$-bimodules $\Hom(W_{\bullet},M)$ and $M\otimes W_{\bullet}$ by a compatible right action of an additional associative algebra $B$. Limiting us to the one vertex case, the DG $\tilde{A}$-bimodules $\Hom(W_{\bullet},M)$ and $M\otimes W_{\bullet}$ are the same as those defined in Subsection \ref{products} when $N=2$. 

Our aim is now to generalize the enriched structure of the DG $\tilde{A}$-bimodule $M\otimes W_{\bullet}$ obtained in the quadratic case~\cite[Section 3]{bt}, to $M\otimes W_{\nu(\bullet)}$ associated with any $N$-homogeneous algebra $A$ -- we leave to the reader the analogous generalization for $\Hom(W_{\nu(\bullet)},M)$. Next, following~\cite[Section 3]{bt}, we will apply this construction to $B=M=A^e$ and then we will show that $K(A)$ has an enriched structure for any 
$N$-homogeneous algebra $A$.

Throughout this section, $A=T(V)/(R)$ denotes an $N$-homogeneous algebra over a field $k$, where $N\geq 2$ and $R$ is a subspace of $V^{\otimes N}$. As usual, we set $A^e=A\otimes A^{op}$ and for any $A$-bimodule $M$, $k$ is assumed to act centrally on $M$. So $M$ can be considered as a left (right) $A^e$-module.

The additional data is the following: $B$ is a unital associative $k$-algebra and the $A$-bimodule $M$ is a right $B$-module such that the actions 
of $k$ induced on $M$ by $A$ and by $B$ are equal. Moreover, the bimodule actions of $A$ and the right action of $B$ on $M$ are assumed to commute, that is, for 
any $m\in M$, one has $(a.m.a').b=a.(m.b).a'$ for $a$, $a'$ in $A$ and $b$ in $B$. Both properties are equivalent to saying that $M$ is an $A^e$-$B$-bimodule. We 
also say that the right $B$-module structure is \emph{compatible} with the $A$-bimodule structure. For example, $M=A^e$ is an $A^e$-$A^e$-bimodule for the product 
of the algebra $A^e$ (here $B=A^e$).

In the sequel, $\catvs$ denotes the category of $k$-vector spaces, $\abimod$ the category of $A$-bimodules, and $\modb$ the category of right-$B$-modules. We have to work with DG modules in abelian categories, following the general framework used by Yekutieli in his recent book on derived categories~\cite{yeku:dercat}. Let us recall from~\cite[Definition 3.8.1]{yeku:dercat} the definition of a DG $\tilde{A}$-bimodule in the abelian category $\modb$. Actually, we extend in an obvious manner this definition to weak DG bimodules over the weak DG algebra $\tilde{A}$. 

\begin{Df} \label{dgbimodin1}
A weak chain (cochain) DG $\tilde{A}$-bimodule $C$ in $\modb$ is a chain (cochain) complex in $\modb$ endowed with a weak DG $\tilde{A}$-bimodule structure such that the bimodule actions of $\tilde{A}$ and the right action of $B$ are compatible. 
\end{Df}

Remark that if the abelian category $\modb$ is replaced by $\catvs$ in this definition, a weak DG $\tilde{A}$-bimodule in $\catvs$ is just a weak DG $\tilde{A}$-bimodule. Remark also that if the weak DG algebra $\tilde{A}$ is replaced by the ground field $k$ (viewed as a trivial DG algebra), we just obtain 
complexes of right $B$-modules.

We start with the weak chain DG $\tilde{A}$-bimodule $M\otimes W_{\nu(\bullet)}$ as recalled at the end of Subsection \ref{products} for any $A$-bimodule $M$. Now $M$ is an $A^e$-$B$-bimodule. We have to define a right action of $B$ on $M\otimes W_{\nu(\bullet)}$. For that, we define a right action of $B$ on each space $M\otimes W_p$, $p\geq 0$, as follows. Fix a basis $(u^{\alpha})_{\alpha}$ of the space $W_p$. Any element $z$ of $M\otimes W_p$ uniquely decomposes as a finitely supported sum $z= \sum_{\alpha} m^{\alpha} \otimes u^{\alpha}$. For any $b\in B$, we set
$$z.b= \sum_{\alpha} (m^{\alpha}.b) \otimes u^{\alpha}.$$
The following lemma shows that the element $z.b$ of $M\otimes W_p$ is well-defined, that is, it does not depend on the choice of the basis $(u^{\alpha})_{\alpha}$.

\Blm \label{lemma1}
For any finite sum $z= \sum_i m_i \otimes w_i$ with $m_i \in M$ and $w_i \in W_p$, one has
\begin{equation} \label{equal}
  \sum_{\alpha} (m^{\alpha}.b) \otimes u^{\alpha}=\sum_{i} (m_i.b)\otimes  w_i.
\end{equation}
\Elm
\Bdm
Fix a basis $(e_{\lambda})_{\lambda}$ of the space $V$. Any $u^{\alpha}$ uniquely decomposes as a finitely supported sum
$$u^{\alpha}= \sum_{\lambda_1 \ldots \lambda_p} c^{\alpha}_{\lambda_1 \ldots \lambda_p} ({e_{\lambda_1} \otimes \cdots \otimes e_{\lambda_p}})$$
where the coefficients are in $k$, so that we have uniquely
\begin{equation} \label{dec1}
  z= \sum_{\alpha, \lambda_1 \ldots \lambda_p} m^{\alpha} c^{\alpha}_{\lambda_1 \ldots \lambda_p} \otimes ({e_{\lambda_1} \otimes \cdots \otimes e_{\lambda_p}}).
\end{equation}
Similarly $w_i$ uniquely decomposes as a finitely supported sum
$$w_i= \sum_{\lambda_1 \ldots \lambda_p} d^i_{\lambda_1 \ldots \lambda_p} ({e_{\lambda_1} \otimes \cdots \otimes e_{\lambda_p}})$$
with coefficients in $k$, so we have uniquely
\begin{equation} \label{dec2}
  z= \sum_{i, \lambda_1 \ldots \lambda_p} m_i d^i_{\lambda_1 \ldots \lambda_p} \otimes ({e_{\lambda_1} \otimes \cdots \otimes e_{\lambda_p}}).
\end{equation}
Comparing (\ref{dec1}) and (\ref{dec2}), we have for any $\lambda_1 \ldots \lambda_p$,
$$\sum_{\alpha} m^{\alpha} c^{\alpha}_{\lambda_1 \ldots \lambda_p}= \sum_{i} m_i d^i_{\lambda_1 \ldots \lambda_p}.$$
If we act with $b$ on this equality and use the compatibility of the actions, we get
$$\sum_{\alpha} (m^{\alpha}.b) c^{\alpha}_{\lambda_1 \ldots \lambda_p}= \sum_{i} (m_i.b) d^i_{\lambda_1 \ldots \lambda_p}.$$
Tensoring by $({e_{\lambda_1} \otimes \cdots \otimes e_{\lambda_p}})$ and summing the so-obtained equalities over the indices $\lambda_1, \ldots \lambda_p$, we obtain formula (\ref{equal}).
\Edm

\medskip

As commonly used in Koszul calculus~\cite{bls:kocal,berger:Ncal,bt}, we write down the element $z$ of $M\otimes W_p$ by 
the compact notation $z=m \otimes x_1 \ldots  x_p$ recalled at the beginning of Subsection \ref{Kbc}. Lemma \ref{lemma1} shows that it makes sense to define 
$z.b$ by the compact notation $z.b=(m.b)\otimes x_1 \ldots  x_p$. 
We systematically use this compact notation in all the sequel. For example, it is clear from this notation that $(z.b).b'=z.(bb')$. 
Notice that the action of $k$ induced by $B$ on $M\otimes W_p$ is equal to the action defined by the $k$-vector space $M\otimes W_p$.

Next, the differential $b_K$ of $M\otimes W_{\nu(\bullet)}$ is right $B$-linear. It suffices to recall~\cite[Subsection 2.2]{berger:Ncal} that, for any $q$-chain $z=m \otimes x_1 \ldots x_{\nu(q)}$, one has 
\begin{equation} \label{defbodd}
b_K(z) =mx_1\otimes x_2 \ldots x_{Nq'+1} - x_{Nq'+1} m \otimes x_1 \ldots x_{Nq'}
\end{equation}
if $q=2q'+1$, and
\begin{equation} \label{defbeven}
b_K(z) = \sum_{0\leq i\leq N-1} x_{i+Nq'-N+2}\ldots x_{Nq'}mx_1\ldots x_i\otimes x_{i+1} \ldots x_{i+Nq'-N+1}
\end{equation}
if $q=2q'$. The compatibility formula $(m.b).a=(m.a).b$ valid for $m \in M$, $a\in A$ and $b \in B$ shows that $b_K(z.b)=b_k(z).b$.

Similarly, for any Koszul $p$-cochain $f:W_{\nu(p)}\rightarrow A$, we easily verify that
$$f \underset{K}{\frown} (z.b)=(f \underset{K}{\frown}z).b\text{ and } (z.b)\underset{K}{\frown} f =(z \underset{K}{\frown}f).b$$
from the defining expressions of $f \underset{K}{\frown}z$ and $z \underset{K}{\frown}f$ recalled in Definition \ref{defcap}. Therefore the bimodule actions of $\tilde{A}$ and the right action of $B$ are compatible. We can conclude by the following.

\Bpo \label{dgalg1}
Let $A=T(V)/(R)$ be an $N$-homogeneous algebra over a field $k$. Fix a unital associative $k$-algebra $B$ and an $A$-bimodule $M$. We assume that $M$ is a right $B$-module compatible with the $A$-bimodule structure. Then the chain complex $M\otimes W_{\nu(\bullet)}$ is a weak DG $\tilde{A}$-bimodule in $\modb$.
\Epo

Let us specialize this result to $B=M=A^e$ as in~\cite[Section 3]{bt}. The left $A^e$-module $A^e$ is identified with $A\stackrel{o}{\otimes} A$, that is, $A\otimes A$ endowed with the outer action $(a\otimes a').(\alpha \otimes \beta)= (a\alpha)\otimes (\beta a')$, where $a$, $a'$, $\alpha$ and $\beta$ are in $A$. Similarly, the right $A^e$-module is identified with $A\stackrel{i}{\otimes} A$ endowed with the inner action $(\alpha \otimes \beta).(a\otimes a')= (\alpha a) \otimes (a' \beta)$. Our aim is to identify the $A$-bimodule complex $K(A)$ with the complex $((A\stackrel{o}{\otimes} A)\otimes W_{\bullet}, b_K)$ endowed with the right action of $A^e$. The following proposition is an $N$-generalization of Proposition 3.5 of~\cite{bt} limited to the one vertex case.

\Bpo
Let $A=T(V)/(R)$ be an $N$-homogeneous algebra over a field $k$.
\begin{enumerate}[\itshape(i)]
\item  For any $q\pgq 0$, the linear map
  $\varphi_q: (A\stackrel{o}{\otimes} A)\otimes W_q \rightarrow A\otimes W_q
  \otimes A$
defined by
  \[\varphi_q ((\alpha \otimes \beta) \otimes x_1 \ldots x_q) = \beta \otimes (x_1 \ldots x_q) \otimes
  \alpha\]
  ia an isomorphism.

\item The direct sum $\Phi$ of the maps $\varphi_{\nu(q)}$ for $q\geq 0$ defines an isomorphism $\Phi$ from the
  complex $((A\stackrel{o}{\otimes} A)\otimes W_{\nu(\bullet)}, b_K)$ to the Koszul complex
  $(K(A),d)$.

\item  The isomorphism $\Phi$ is right $A^e$-linear.
\end{enumerate}
\Epo

\Bdm
The assertion \textit{(i)} is clear. Let us show that $\Phi$ is a morphism of complexes. Take $z=(\alpha \otimes \beta) \otimes x_1 \ldots x_{\nu(q)}$. If $q=2q'+1$, we deduce from
$$b_K(z) =(\alpha \otimes \beta x_1)\otimes x_2 \ldots x_{Nq'+1} - (x_{Nq'+1}\alpha\otimes \beta) \otimes x_1 \ldots x_{Nq'},$$
the following
$$\Phi(b_K(z))=\beta x_1 \otimes x_2 \ldots x_{Nq'+1} \otimes \alpha - \beta  \otimes x_1 \ldots x_{Nq'}\otimes x_{Nq'+1} \alpha$$
whose right-hand side is equal to $d(\beta \otimes x_1 \ldots x_{Nq'+1} \otimes \alpha)$.

Now if $q=2q'$, from
$$b_K(z) = \sum_{0\leq i\leq N-1} (x_{i+Nq'-N+2}\ldots x_{Nq'}\alpha \otimes \beta x_1\ldots x_i) \otimes x_{i+1} \ldots x_{i+Nq'-N+1}$$
we get
$$\Phi(b_K(z))= \sum_{0\leq i\leq N-1} \beta x_1\ldots x_i \otimes x_{i+1} \ldots x_{i+Nq'-N+1} \otimes x_{i+Nq'-N+2}\ldots x_{Nq'}\alpha$$
whose right-hand side is equal to $d(\beta \otimes x_1 \ldots x_{Nq'} \otimes \alpha)$.

Let us prove \textit{(iii)}. Here the $A$-bimodule $A\otimes W_q \otimes A$ is seen as a right $A^e$-module.
For $z=(\alpha \otimes \beta) \otimes x_1 \ldots x_q$ and $a$, $a'$ in $A$, we have 
\begin{align*}
\varphi_q(z.(a\otimes a'))&=\varphi_q((\alpha a \otimes a' \beta) \otimes x_1 \ldots x_q)\\&=a' \beta \otimes (x_1  \ldots x_q) \otimes \alpha a\\&= \varphi_q(z).(a\otimes a'),
\end{align*} therefore $\varphi_q$ is right $A^e$-linear for any $q\geq 0$. In particular, $\Phi_q=\varphi_{\nu(q)}$ is right $A^e$-linear for any $q\geq 0$. 
\Edm

\medskip

By Proposition \ref{dgalg1}, the chain complex $(A\stackrel{o}{\otimes} A)\otimes W_{\bullet}$ is a weak DG $\tilde{A}$-bimodule in $\abimod$. We transport this structure via the chain complex isomorphism $\Phi$ and we obtain.

\begin{Po} \label{dgalg2}
Let $A=T(V)/(R)$ be an $N$-homogeneous algebra over a field $k$. Then the Koszul bimodule complex $K(A)$ is a weak DG $\tilde{A}$-bimodule in $\abimod$.
\end{Po}

Let us give explicitly the underlying weak $\tilde{A}$-bimodule structure of the weak DG $\tilde{A}$-bimodule $K(A)$. Consider $z=(\alpha \otimes \beta) \otimes x_1 \ldots x_{\nu(q)}$ in $(A\stackrel{o}{\otimes} A)\otimes W_{\nu(q)}$ and $f$ in $\Hom (W_{\nu(p)},A)$. Set $z'= \Phi_q(z)= \beta \otimes x_1 \ldots  x_{\nu(q)} \otimes \alpha$.

If $p$ and $q-p$ are not both odd, we draw from the formulas giving  $f\underset{K}{\frown} z$ and $z\underset{K}{\frown} f$ in Definition \ref{defcap} that the actions of $f$ on $K(A)$ are defined by
\begin{equation} 
f \underset{K}{\frown} z' = (-1)^{(q-p)p} \beta \otimes x_1 \ldots  x_{\nu(q-p)} \otimes f(x_{\nu(q-p)+1} \ldots x_{\nu(q)})\alpha,
\end{equation}
\begin{equation} 
z' \underset{K}{\frown} f = (-1)^{pq} \beta f(x_1 \ldots x_{\nu(p)})\otimes x_{\nu(p)+1} \ldots  x_{\nu(q)} \otimes \alpha.
\end{equation}

If $p=2p'+1$ and $q=2q'$, we use again the formulas giving  $f\underset{K}{\frown} z$ and $z\underset{K}{\frown} f$ in Definition \ref{defcap}, and we easily deduce
\begin{eqnarray*}
  f\underset{K}{\frown} z' =  -\sum_{0\leq i+j \leq N-2} (\beta x_1 \ldots x_i)\otimes (x_{i+1} \ldots  x_{i+Nq'-Np'-N+1})\\ \otimes (x_{Nq'-Np'-N+i+2}\ldots x_{Nq'-Np'-j-1} f(x_{Nq'-Np'-j} \ldots x_{Nq'-j}) x_{Nq'-j+1}\ldots x_{Nq'}\alpha),
\end{eqnarray*}
\begin{eqnarray*}
  z'\underset{K}{\frown} f =  \sum_{0\leq i+j \leq N-2} (\beta x_1 \ldots x_i f(x_{i+1} \ldots x_{Np'+i+1}) x_{Np'+i+2}\ldots x_{Np'+N-j-1})\\ \otimes (x_{Np'+N-j}\ldots x_{Nq'-j})\otimes (x_{Nq'-j+1}\ldots x_{Nq'}\alpha).
\end{eqnarray*}

\setcounter{equation}{0}

\section{Towards a proof of the graded commutativity} \label{towards}

\subsection{Introduction}

We would like to obtain a proof of the graded commutativity of the Koszul cup product for any $N$-homogeneous algebra $A$, that is, a proof of property (i) in Corollary \ref{cupcapcons}. Our idea is to divide such a proof in two steps as follows.

1) First Step: use the fact that the Koszul cohomology is isomorphic to a Hochschild hypercohomology. This fact was obtained when $N=2$ in~\cite[Subsection 2.3]{bls:kocal}, and extended without detailed proof to any $N$ in~\cite[Subsection 2.2]{berger:Ncal}.

2) Second Step: show that the isomorphism of the First Step sends the Koszul cup product to a Hochschild cup product. Then we could conclude by using the graded commutativity of the Hochschild cup product~\cite{gerst:cohom}.

The First Step is stated in Proposition \ref{isohyper} below, including a detailed proof. The Second Step lies on the general strategy presented in the previous section. The enriched structures will be explicitly described for any $N$-homogeneous algebra $A$ in Subsection \ref{enriched} below. An extra hypothesis (H) is needed even if $N=2$ in order to make use of derived categories (Theorem \ref{derivedrho}). 
Moreover we obtain a product on the image of the isomorphism that we would still need to relate to a usual Hochschild cup product.

\subsection{Koszul cohomology and Hochschild hypercohomology}

\Bpo \label{isohyper}
Let $A= T(V)/(R)$ be an $N$-homogeneous algebra over a field $k$. For any $A$-bimodule $M$, the Koszul cohomology $\HK^{\bullet}(A,M)$ is isomorphic, as a graded vector space, to the Hochschild hypercohomology $\mathbb{H}\mathbb{H}^{\bullet}(A, \Hom _A (K(A),M))$.
\Epo
\Bdm
For the definition of the derived category $\mathcal{D} (\mathcal{C})$ of an abelian category $\mathcal{C}$, we refer to~\cite[Chapter 10]{weib:homo}. We will make precise below the boundedness conditions which we will use. Denote by $\catvs$ the category of $k$-vector spaces and by $\abimod$ the category of $A$-bimodules. 

Our proof is based on an isomorphism of functors depending on a fixed $A$-bimodule $M$ that we will consider as left $A^e$-module. In general, unless 
the contrary is explicitly stated, $A$-bimodules will always be identified to left $A^e$-modules.
Let $F: \abimod \rightarrow \catvs$ be the contravariant functor $F: P \mapsto \Hom _{A^e}(P,M)$.
 Let $G: \abimod \rightarrow \abimod$ be the contravariant functor $G: P\mapsto \Hom_{A} (P,M)$ where $\Hom_{A} (P,M)$ denotes the space of left $A$-module morphisms $u:P\rightarrow M$, 
endowed with the following $A$-bimodule structure
\begin{equation} \label{actions}
(a.u.a')(x)= u(xa)a', \hbox{ for }  a, a'\in A, \ x\in P.
\end{equation}
Finally, let $H: \abimod \rightarrow \catvs$ be the covariant functor $H: P \mapsto \Hom _{A^e}(A,P)$.

Denote by $\mathcal{C_+}(\abimod)$ the category of chain complexes $C=(C_n)$ of $A$-bimodules which are bounded below, that is, such that $C_n=0$ for all $n<<0$. Denote by $\mathcal{C^+}(\abimod)$ and $\mathcal{C^+}(\catvs)$ the categories of the cochain complexes $C=(C^n)$ of $A$-bimodules 
--respectively, vector spaces-- which are bounded below, that is, such that $C^n=0$ for all $n<<0$. Then we have natural functors 
that we deduce from $F$, $G$, $H$ and we still denote as before 
$F:\mathcal{C_+} (\abimod) \rightarrow \mathcal{C^+} (\catvs)$, $G: \mathcal{C_+} (\abimod) \rightarrow \mathcal{C^+} (\abimod)$ and $H:\mathcal{C^+} (\abimod) \rightarrow \mathcal{C^+} (\catvs)$.

Following~\cite[Chapter 10]{weib:homo}, we consider the right derived functors $RF:\mathcal{D_+} (\abimod) \rightarrow \mathcal{D^+} (\catvs)$, $RG:\mathcal{D_+} (\abimod) \rightarrow \mathcal{D^+} (\abimod)$ and $RH:\mathcal{D^+} (\abimod) \rightarrow \mathcal{D^+} (\catvs)$.

We define a natural transformation of functors $\rho_M: F\to H\circ G$ as follows. For any $A$-bimodule $P$, we first 
define the map 
$$\rho_M^P: \Hom_{A^e}(P, M) \rightarrow \Hom_{A^e} (A,Hom_A(P,M))$$
by $\rho_M^P(u)(a)(x)=u(xa)$ for any $A$-bimodule morphism $u:P\rightarrow M$, $a$ in $A$ and $x$ in $P$. 
It is easy to check that $\rho_M^P(u)$ belongs to $\Hom_{A^e} (A,Hom_A(P,M))$. Since $\rho_M^P(u)(1)=u$, it follows that $\rho_M^P$ is linear. For any $v$ in 
$\Hom_{A^e} (A,\Hom_A(P,M))$, define $\rho_M'^P(v)=v(1)$. Then $\rho_M'^P(v)$ belongs to $Hom_{A^e}(P, M)$ and the map $\rho_M'^P$ is linear. Consequently, 
the maps $\rho_M^P$ and $\rho_M'^P$ are inverse to each other. So $\rho_M^P$ is a linear isomorphism.

Moreover $\rho_M^P$ is functorial in $P$, defining thus a natural isomorphism 
$\rho_M: F \cong H \circ G$ from $\mathcal{C_+}(\abimod)$ to $\mathcal{C^+}(\catvs)$. Therefore, $R(\rho_M)$ is an isomorphism of functors $RF \cong RH \circ RG$ from $\mathcal{D_+} (\abimod)$ to $\mathcal{D^+} (\catvs)$.  In particular, for any bounded below chain complex $C$ of projective $A$-bimodules, one has an isomorphism 
\begin{equation} \label{isocostand}
R(\rho_M)(C): \RHom_{A^e}(C, M) \cong \RHom_{A^e}(A, \RHom_{A}(C, M))
\end{equation}
in $\mathcal{D^+} (\catvs)$. Applying this isomorphism to $C=K(A)$, passing to cohomology and using the definition of the hypercohomology~\cite[Chapter 10]{weib:homo}, we obtain a graded linear isomorphism
\begin{equation} \label{hyperhcoh}
\HK^{\bullet}(A,M) \cong \mathbb{H}\mathbb{H}^{\bullet}(A, \Hom_{A}(K(A), M)).
\end{equation}
Remark that, if $A$ is $N$-Koszul, then $K(A) \cong A$ in $\mathcal{D_+} (\abimod)$ so that we recover $\HK^{\bullet}(A,M) \cong \HH^{\bullet}(A,M))$, see~\cite[Subsection 2.2]{berger:Ncal}.
\Edm

\subsection{Koszul cohomology and Hochschild hypercohomology, with enriched structures} \label{enriched}

Let us present precisely our strategy for proving property (i) in Corollary \ref{cupcapcons} for any $N$-homogeneous algebra $A= T(V)/(R)$. We have to work with weak DG $\tilde{A}$-bimodules in the abelian category $\abimod$. As seen in Proposition \ref{dgalg2}, $K(A)$ is such a weak DG $\tilde{A}$-bimodule in $\abimod$. Recall Definition \ref{dgbimodin1} in the particular case $B=A^e$.

\begin{Df} \label{dgbimodin2}
A weak chain (cochain) DG $\tilde{A}$-bimodule $C$ in $\abimod$ is a chain (cochain) complex in $\abimod$ endowed with a weak DG $\tilde{A}$-bimodule structure such that the bimodule actions of $\tilde{A}$ and $A$ are compatible. 
\end{Df}

Denote by $\mathcal{C}_+^w(\tilde{A},\abimod)$ the category of weak bounded below chain DG $\tilde{A}$-bimodules in $\abimod$. Denote by $\mathcal{C}_w^+(\tilde{A},\abimod)$ and $\mathcal{C}_w^+(\tilde{A},\catvs)$ the category of weak bounded below cochain DG $\tilde{A}$-bimodules in $\abimod$ and $\catvs$ respectively. Notice that $\mathcal{C}_w^+(\tilde{A},\catvs)$ coincides with the category $\mathcal{C}_w^+(\tilde{A})$ of weak bounded below cochain DG $\tilde{A}$-bimodules.

In case $N=2$, $K(A)$ is a DG $\tilde{A}$-bimodule in $\abimod$, so that in the above notations we can drop the subscript $w$ and we have 
that $\mathcal{C_+}(\tilde{A},\abimod)$ is the category of bounded below chain DG $\tilde{A}$-bimodules in $\abimod$, while $\mathcal{C^+}(\tilde{A},\abimod)$ and $\mathcal{C^+}(\tilde{A},\catvs)=\mathcal{C^+}(\tilde{A})$ are the categories of bounded below cochain DG $\tilde{A}$-bimodules in $\abimod$ and $\catvs$ respectively. 

We begin with the general case $N\geq 2$. Fix an $A$-bimodule $M$ and a weak chain DG $\tilde{A}$-bimodule $C$ in $\abimod$. It is easy to check that $\Hom _{A^e} (C, M)$, where $\Hom_{A^e}(C,M)^q:=\Hom _{A^e} (C_q,M)$, is a weak cochain DG $\tilde{A}$-bimodule in $\catvs$ for the actions
\begin{equation} \label{actions2}
  (f.u)(x)=(-1)^p u(x.f), \ \ (u.f)(x)=u(f.x),
\end{equation}
where $f:A\otimes_k W_p \otimes_k A \rightarrow A$, $u: C_q \rightarrow M$ are morphisms of $A$-bimodules, and $x \in C_{p+q}$. Note that $x.f$ and $f.x$ are in $C_q$ by the graded actions of $\tilde{A}$ on $C$. Recall also that if $d:C_{q+1}\rightarrow C_q$ is the differential of $C$, the differential
$$\Hom_{A^e}(d,M): \Hom _{A^e} (C_q,M) \rightarrow \Hom _{A^e} (C_{q+1},M)$$
of $\Hom _{A^e} (C,M)$ is defined by $\Hom_{A^e}(d,M)(u)=-(-1)^q u\circ d$ for $u \in \Hom _{A^e} (C_q,M)$. So we have defined a functor
$$\Hom_{A^e}(-,M) :\mathcal{C}_+^w (\tilde{A},\abimod) \rightarrow \mathcal{C}^+_w (\tilde{A},\catvs).$$
Note that if $C=K(A)$, then $x.f= x\underset{K}{\frown}f$ and $f.x=f \underset{K}{\frown}x$ are explicitly expressed just after Proposition \ref{dgalg2}, so that we recover the Koszul cup actions, that is, $f.u=f\underset{K}{\smile}u$ and $u.f=u\underset{K}{\smile}f$. Remark also that, if $N=2$ and $C$ is a chain DG $\tilde{A}$-bimodule in $\abimod$, then $\Hom _{A^e} (C, M)$ is a cochain DG $\tilde{A}$-bimodule in $\catvs$, as seen in~\cite[Subsection 5.3]{bt}.

Keeping the same notation but now $u: C_q \rightarrow M$ is only left $A$-linear, we assert that $\Hom _{A} (C, M)$, where $\Hom_{A}(C,M)^q:=\Hom _{A} (C_q,M)$, is a weak cochain DG $\tilde{A}$-bimodule in $\abimod$. Actually, we just explain the structures involved in Definition \ref{dgbimodin2} and we leave the verifications of their properties to the reader. Firstly the differential 
$$\Hom_{A}(d,M): \Hom _{A} (C_q,M) \rightarrow \Hom _{A} (C_{q+1},M)$$
of $\Hom _{A} (C,M)$ is defined by $\Hom_{A}(d,M)(u)=-(-1)^q u\circ d$ for $u \in \Hom _{A} (C_q,M)$. It is a differential of $A$-bimodules for the actions defined in (\ref{actions}). Secondly, the actions of $\tilde{A}$ being defined by the same formulas (\ref{actions2}), it is easy to verify that $f.u$ and $u.f$ are left $A$-linear. Finally it is routine to check that the actions of $A$ and $\tilde{A}$ on $\Hom _{A} (C,M)$ are compatible, and that the following formulas hold
\begin{align}
  \label{derivleft}
   \Hom_{A}(d,M)(f.u) &= b_K(f).u + (-1)^p f.\Hom_{A}(d,M)(u),\\
  \label{derivright}
  \Hom_{A}(d,M)(u.f) &= \Hom_{A}(d,M)(u).f + (-1)^q u.b_K(f).
  \end{align}
In conclusion, we have defined a functor
$$\Hom_{A}(-,M) :\mathcal{C}_+^w (\tilde{A},\abimod) \rightarrow \mathcal{C}^+_w (\tilde{A},\abimod).$$
Note that if $C=K(A)$ and $u$ is only left $A$-linear, we have only $f.u=f\underset{K}{\smile}u$. Remark again that, if $N=2$ and $C$ is a chain DG $\tilde{A}$-bimodule in $\abimod$, we can remove ``weak'' everywhere as well as the subscript $w$ (we leave the verifications to the reader).

Similarly, for any weak cochain DG $\tilde{A}$-bimodule $(C,d)$ in $\abimod$, $\Hom_{A^e}(A, C)$, where $\Hom_{A^e}(A,C)^q:= \Hom_{A^e}(A,C^q)$, is a weak cochain DG $\tilde{A}$-bimodule in $\catvs$. The differential $\Hom_{A^e}(A,d)$ is defined by
$$\Hom_{A^e}(A,d)(v)=d\circ v$$
for $v:A\rightarrow C$. The actions of $f:A\otimes_k W_p \otimes_k A \rightarrow A$ on $v:A\rightarrow C^{p+q}$ are defined by
\begin{equation} \label{actions3}
  (f.v)(a)=f.v(a), \ \ (v.f)(a)=v(a).f,
\end{equation}
for any $a \in A$. It is straightforward to show that we have so obtained a functor
$$\Hom_{A^e}(A,-) :\mathcal{C}^+_w (\tilde{A},\abimod) \rightarrow \mathcal{C}^+_w (\tilde{A},\catvs).$$
Here again, if $N=2$ and $C$ is a cochain DG $\tilde{A}$-bimodule in $\abimod$, we can remove ``weak'' and the subscript $w$ everywhere.

\Bpo \label{enrichedrho}
Let $A= T(V)/(R)$ be an $N$-homogeneous algebra over a field $k$ and $M$ be an $A$-bimodule. The isomorphism of functors $\rho_M: F\cong H\circ G$ defined in the proof of Proposition \ref{isohyper} has an enriched version
$$\rho_M: \Hom_{A^e}(-,M) \cong \Hom_{A^e}(A,-) \circ \Hom_{A}(-,M)$$
which is an isomorphism of functors from $\mathcal{C}_+^w (\tilde{A},\abimod)$ to $\mathcal{C}^+_w (\tilde{A})$. If $N=2$, one has an enriched version
$$\rho_M: \Hom_{A^e}(-,M) \cong \Hom_{A^e}(A,-) \circ \Hom_{A}(-,M)$$
which is an isomorphism of functors from $\mathcal{C}_+ (\tilde{A},\abimod)$ to $\mathcal{C}^+(\tilde{A})$.

\Epo
\Bdm
Fix a weak chain DG $\tilde{A}$-bimodule $(C,d)$ in $\abimod$. It is easy to check that 
$$\rho_M^C: \Hom_{A^e}(C,M) \cong \Hom_{A^e}(A,\Hom_{A}(C,M))$$
is a morphism of complexes. It remains to verify that $\rho_M^C$ is a morphism of weak $\tilde{A}$-bimodules. We use the structures of weak $\tilde{A}$-bimodules introduced above. Consider the following morphisms of $A$-bimodules: $f:A\otimes_k W_p \otimes_k A \rightarrow A$, $f':A\otimes_k W_{p'} \otimes_k A \rightarrow A$ and $u: C_q \rightarrow M$. Then $f.u.f' \in \Hom_{A^e}(C_{p+p'+q},M)$ and we take $a \in A$ and $x \in C_{p+p'+q}$.

On one hand, one has $\rho_M^C(f.u.f')(a)(x)=(-1)^p u(f'.(xa).f)$. On the other hand, one has $(f.\rho_M^C(u).f')(a)=f.\rho_M^C(u)(a).f'$ and
$$(f.\rho_M^C(u)(a).f')(x)=(-1)^p\rho_M^C(u)(a)(f'.x.f)= (-1)^p u((f'.x.f)a).$$
But the bimodule actions of $\tilde{A}$ and of $A$ on $C$ are compatible. In particular, $f'.(xa).f= (f'.x.f)a$, so that we conclude that $\rho_M^C(f.u.f')=f.\rho_M^C(u).f'$.

If $N=2$, we can remove ``weak'' and the subscript $w$ everywhere.
\Edm
\\

We do not know whether the corresponding derived categories make sense if $N>2$. However, if $N=2$, we can remove ``weak'' and the 
subscript $w$, and so we recover the framework of Yekutieli~\cite{yeku:dercat}. Consequently, we assume that $N=2$ throughout the remainder of this section. 
In this situation, the derived categories $\mathcal{D_+} (\tilde{A},\abimod)$, $\mathcal{D^+} (\tilde{A},\abimod)$ and $\mathcal{D^+} (\tilde{A},\catvs)=\mathcal{D^+} (\tilde{A})$ are defined in~\cite{yeku:dercat}. We are now ready to introduce the following hypothesis (H).
\\

{\em Hypothesis} (H): The right derived functors $\RHom_{A^e}(-,M):\mathcal{D_+} (\tilde{A},\abimod) \rightarrow \mathcal{D^+} (\tilde{A})$ and $\RHom_{A}(-,M):\mathcal{D_+} (\tilde{A},\abimod) \rightarrow \mathcal{D^+} (\tilde{A},\abimod)$ exist for any $A$-bimodule $M$, as well as the right derived functor 
$\RHom_{A^e}(A,-):\mathcal{D^+} (\tilde{A},\abimod) \rightarrow \mathcal{D^+} (\tilde{A})$. 
\\

In his book, Yekutieli shows that if suitable resolutions exist, then the functors between categories of DG modules in an abelian category can be derived~\cite[Theorem 10.1.20, Theorem 10.2.15, Theorem 10.4.8, Theorem 10.4.9]{yeku:dercat}. We do not know whether these general results by Yekutieli can be applied, in the present context, to the existence of injective resolutions for proving the assertion (H). However, from Proposition \ref{enrichedrho} and from general properties of abstract derived functors~\cite[Subsection 8.3]{yeku:dercat}, we obtain the following.

\Bte \label{derivedrho}
Let $A= T(V)/(R)$ be a quadratic algebra over a field $k$ and $M$ be an $A$-bimodule. Under the hypothesis (H), the isomorphism of functors $\rho_M: \Hom_{A^e}(-,M) \cong \Hom_{A^e}(A,-) \circ \Hom_{A}(-,M)$
from $\mathcal{C_+} (\tilde{A},\abimod)$ to $\mathcal{C^+} (\tilde{A})$ induces an isomorphism of functors
$$R(\rho_M): \RHom_{A^e}(-,M) \cong \RHom_{A^e}(A,-) \circ \RHom_{A}(-,M)$$
from $\mathcal{D_+} (\tilde{A},\abimod)$ to $\mathcal{D^+} (\tilde{A})$.
\Ete

In particular, for any DG $\tilde{A}$-bimodule $C$ in $\abimod$ which is a bounded below chain complex of projective $A$-bimodules, we obtain an isomorphism 
\begin{equation} \label{isocostand2}
\RHom_{A^e}(C, M) \cong \RHom_{A^e}(A, \RHom_{A}(C, M))
\end{equation}
in $\mathcal{D^+} (\tilde{A})$. Applying this isomorphism to $C=K(A)$ and passing to cohomology, we get that the graded linear isomorphism (\ref{hyperhcoh}), that is, 
$$\HK^{\bullet}(A,M) \cong \mathbb{H}\mathbb{H}^{\bullet}(A, \Hom _A (K(A),M))$$
is now an isomorphism of graded $\HK^{\bullet}(A)$-bimodules. 

Therefore, in order to prove that the graded $\HK^{\bullet}(A)$-bimodule $\HK^{\bullet}(A,M)$ is commutative, it suffices to prove that the graded $\HK^{\bullet}(A)$-bimodule $\mathbb{H}\mathbb{H}^{\bullet}(A, \Hom _A (K(A),M))$ is commutative. Unfortunately, we have not succeeded to relate the so-obtained graded action of $\HK^{\bullet}(A)$ on $\mathbb{H}\mathbb{H}^{\bullet}(A, \Hom _A (K(A),M))$ to the usual Hochschild cup product which is known to be graded commutative.

\vspace{0.5 cm} \textsf{Roland Berger: Univ Lyon, UJM-Saint-\'Etienne, CNRS UMR 5208, Institut Camille Jordan, F-42023, Saint-\'Etienne, France}

\emph{roland.berger@univ-st-etienne.fr}\\

\textsf{Andrea Solotar: IMAS and Dto de Matem\'{a}tica, Facultad de Ciencias Exactas y Naturales,
Universidad de Buenos Aires, Ciudad Universitaria, Pabell\`{o}n 1,
(1428) Buenos Aires, Argentina}

\emph{asolotar@dm.uba.ar}

\end{document}